\documentclass[preprints,article,submit,moreauthors,pdftex]{Definitions/mdpi}

\firstpage{1}
\pubvolume{1}
\issuenum{1}
\articlenumber{0}
\pubyear{2022}
\copyrightyear{2022}
\datereceived{}
\dateaccepted{}
\datepublished{}
\hreflink{https://doi.org/} 
\pdfoutput=1

\usepackage{amssymb}

\usepackage[english]{babel}

\def\dtp#1{\mathop {#1}\limits_{+\tau}}
\def\dtm#1{\mathop {#1}\limits_{-\tau}}
\def\dsp#1{\mathop {#1}\limits_{+s}}
\def\dsm#1{\mathop {#1}\limits_{-s}}
\newcommand{\unknownCL}{{Unknown}}

\preto{\abstractkeywords}{\nolinenumbers}

\Title{Invariant Finite-Difference Schemes With Conservation Laws Preservation For One-Dimensional~MHD Equations}

\TitleCitation{Invariant Schemes For MHD}

\Author{Evgeniy Kaptsov $^{1,\textmd{*}}$\orcidB{}, Vladimir Dorodnitsyn $^{2}$\orcidA{}}

\AuthorNames{Evgeniy Kaptsov, Vladimir Dorodnitsyn}

\AuthorCitation{Kaptsov, E.; Dorodnitsyn, V.}

\address{%
$^{1}$ \quad School of Mathematics, Institute of Science,
Suranaree University of Technology, Nakhon Ratchasima, 30000, Thailand; evgkaptsov@math.sut.ac.th, evgkaptsov@gmail.com\\
$^{2}$ \quad Keldysh Institute of Applied Mathematics,
Russian Academy of Science, Miusskaya Pl. 4, Moscow, 125047, Russia; Dorodnitsyn@keldysh.ru, dorod2007@gmail.com}

\corres{Correspondence: evgkaptsov@gmail.com}

\abstract{Invariant finite-difference schemes are considered for
one-dimensional magnetohydrodynamics~(MHD) equations in
mass~Lagrangian coordinates for the cases of finite and infinite
conductivity. The construction of these schemes make use of previously obtained
results of the group classification of MHD equations~\cite{bk:DorKozMelKap_PlainFlows_2021}. 
On the basis of the classical Samarskiy--Popov scheme, new schemes
are constructed for the case of finite conductivity. These schemes
admit all symmetries of the original differential model and have
difference analogues of all of its local differential conservation
laws. 
New, previously unknown, conservation laws are found using symmetries
and direct calculations. 
In the case of infinite conductivity, conservative invariant schemes are
constructed as well. For isentropic flows of a polytropic
gas the proposed schemes possess the conservation law of energy and
preserve entropy on two time layers. This is achieved by means of
specially selected approximations for the equation of state of a
polytropic gas. Also, invariant difference schemes with additional
conservation laws are proposed. 
A new scheme for the case of finite conductivity is tested
numerically for various boundary conditions which shows accurate
preservation of difference conservation laws.} 

\keyword{classical symmetries; conservation law; numerical scheme}

\begin{document}

\section{Introduction}

Magnetic hydrodynamics equations describe the flows of electrically conductive fluids
such as plasma, liquid metals, and electrolytes
and are widely used in modeling processes
in various fields from engineering to geophysics and astrophysics.

In the present publication, we restrict ourselves to considering plane one-dimensional MHD flows
under the assumption that the medium is inviscid and thermally non-conducting.
A group classification of the MHD equations under the above conditions was carried out recently in~\cite{bk:DorKozMelKap_PlainFlows_2021} (for some particular results see also~\cite{bk:Gridnev1968,bk:DorMHDpreprint1976,art:Rogers1969,bk:HandbookLie_v2}).
The group classification splits into four essentially different cases
according to whether the conductivity of the medium is finite or infinite,
and the longitudinal component of the magnetic field vector is zero or a non-zero constant.

\medskip

The MHD equations are nonlinear, so that even in the one-dimensional case only
their particular solutions are known~\cite{bk:Oliveri_F_2005,bk:Picard_P_Y_2008,bk:Golovin_2009,bk:Golovin_2011,bk:Golovin_2019}.
Therefore, numerical modeling in magnetohydrodynamics is of great practical interest.
There are many approaches to numerically modeling MHD equations, including finite-difference, finite
element and finite volume methods~(see, e.g.~\cite{Toro,
bk:SamarskyPopov[1970],
bk:SamarskyPopov_book[1992],
bk:FalleKomissarovJoarderMHD,
bk:POWELL1999284,bk:YAKOVLEV201380,
bk:YANG2017561,bk:HIRABAYASHI2016851,bk:DongsuMHD}).
Further we consider finite-difference schemes taking as a starting point the classical Samarsky--Popov schemes~\cite{bk:SamarskyPopov[1970],bk:SamarskyPopov_book[1992]} for the MHD equations for the case of finite conductivity.
The main properties of the considered schemes are invariance, i.e. preservation of the symmetries of
the original differential equations, and the presence of difference analogues of local differential conservation laws.
It is known that there is a connection between the invariance of equations and the presence of conservation laws~\cite{bk:Noether1918,bk:Ibragimov1985,bk:Olver,bk:BlumanAnco}.

\medskip

Invariant schemes have been studied for a long time~\cite{Dor_1,[Dor_4],Maeda1,Maeda2,bk:Dorodnitsyn[2011]}, and over the past decades,
significant progress has been made in the development of methods for their construction and integration.
For schemes for ordinary differential equations with Lagranagian or Hamiltonian functions,
a number of methods~\cite{bk:DorodKozlovWint[2004],bk:DorodnitsynKozlovWint[2003],bk:Dorod_Hamilt[2009],bk:Dorod_Hamilt[2010]} have been developed that make it possible to decrease the order or even integrate the schemes.
A method based on the Lagrangian identity has also been developed for the case when the equations do not admit a variational formulation~\cite{bk:DorodKozlovWintKaptsov[2014],bk:DorodKozlovWintKaptsov[2015]}.

For partial difference schemes, the main methods used are the method
of differential
invariants~\cite{bk:Dorodnitsyn[2011],[Pavel],[Pavel2]} and the
difference analogue of the direct
method~\cite{bk:Dorodnitsyn[2011],bk:ChevDorKap2020}. Using these
methods, the authors have constructed invariant schemes for various
shallow water
models~\cite{dorodnitsyn2019shallow,bk:DorKapMelGN2020,DorKapSWJMP2021,bk:DorKapMel_ModSW_2021}.
Also, some previously known schemes have been investigated from a
group analysis point of view. In particular, symmetries and
conservation laws of the Samarskiy--Popov schemes for the
one-dimensional gas dynamics equations of a polytropic gas have been
investigated in~\cite{DORODNITSYN2019201,KOZLOV2019,bk:DorKozMelGasdynBookChapter2021}. Based on the
results of the group
classification~\cite{bk:DorKozMelKap_PlainFlows_2021} and
Samarskiy--Popov schemes for the MHD equations, we further construct
invariant finite-difference schemes possessing conservation laws.
The set and number of conservation laws  depend on the conductivity,
the form of the magnetic field vector and the equation of state of
the medium.

\bigskip

This paper is organized as follows.
In Section~\ref{sec:basics} the simplest version of the finite-conductivity MHD equations in mass Lagrangian coordinates in case of one-dimensional plane flows is considered.
Electric and magnetic fields are represented by one-component vectors, which greatly simplifies the form of the equations.
This was the main case considered in Samarsky and Popov's publications~\cite{bk:SamarskyPopov[1970],bk:SamarskyPopov_book[1992]}.
The section also provides basic notation and definitions.
Then, symmetries and conservation laws of the Samarsky-Popov scheme for the MHD equations are investigated.
In addition to the previously known conservation laws, the center-of-mass conservation law is given,
as well as new conservation laws for the specific conductivity function,
obtained on the basis of the group classification~\cite{bk:DorKozMelKap_PlainFlows_2021}.

Section~\ref{sec:general} is devoted to various generalizations of
the scheme of Section~\ref{sec:basics}. In
Section~\ref{sec:genFiniteCond} the scheme for arbitrary electric
and magnetic fields is considered. Its symmetries are investigated
and conservation laws are given. The case of infinite conductivity
is considered in Section~\ref{sec:genInfiniteCond}. It is shown that
in this case the Samarsky-Popov scheme requires some additional
modifications in order to possess the conservation law of angular
momentum. In the case of a polytropic gas, it turns out to be
possible to preserve not only energy, but also entropy along
pathlines. This can be done using a specially selected equation of
state for a polytropic gas.
At the end of the section, an example of
an invariant scheme is given that does not possess a conservation
law of energy, but preserves entropy and has additional conservation
laws in the case of isentropic flows.
In Section~\ref{sec:experiments}, one of the invariant schemes for the case of finite conductivity is numerically implemented for the example of plasma bunch deceleration by crossed electromagnetic fields.
The results are discussed in
the~Conclusion.

\section{Conservative schemes for MHD equations with finite conductivity}
\label{sec:basics}

Problems of continuum mechanics and plasma physics are often
considered in mass Lagrangian coordinates
~\cite{bk:SamarskyPopov_book[1992],bk:YanenkRojd[1968]} since for
them the formulation of boundary conditions is greatly simplified.
In particular, the conservative Samarsky--Popov schemes for the
equations of gas dynamics and magnetohydrodynamics have been
constructed in mass Lagrangian coordinates.

\medskip

In mass Lagrangian coordinates the MHD equations, describing the
plane one-dimensional MHD flows, are~\cite{bk:SamarskyPopov_book[1992],bk:DorKozMelKap_PlainFlows_2021}
\begin{subequations} \label{eq:sys1ext}
\begin{equation}
\left(\frac{1}{\rho}\right)_t = u_s,
\end{equation}
\begin{equation}
u_t = -\left( p + \frac{(H^y)^2 + (H^z)^2}{2} \right)_s,
\end{equation}
\begin{equation}
v_t = H^0 H^y_{{s}},
\qquad
y_t = v,
\end{equation}
\begin{equation}
w_t = H^0 H^z_{{s}},
\qquad
z_t = w,
\end{equation}
\begin{equation}
\left(\frac{H^y}{\rho}\right)_t = (H^0 v + E^z)_s,
\end{equation}
\begin{equation}
\left(\frac{H^z}{\rho}\right)_t = (H^0 w - E^y)_s,
\end{equation}
\begin{equation}
\sigma E^y = -\rho H^z_{{s}},
\qquad
\sigma E^z = \rho H^y_{{s}},
\end{equation}
\begin{equation}
\label{eq:energyVec}
\varepsilon_t = -p u_s + \frac{1}{\rho}(\mathbf{i} \cdot \mathbf{E}),
\end{equation}
\begin{equation}
x_t = u,
\qquad
x_s = 1/\rho,
\end{equation}
\end{subequations}
where
$t$ is time,
$s$ is mass Lagrangian coordinate,
$x$ is Eulerian coordinate,
$\rho$ is density,
$p$ is pressure,
$\varepsilon$ is internal energy,
$\mathbf{u}=(u,v,w)$ is the velocity of a particle,
$\mathbf{E}=(E^x,E^y,E^z)$ is the electric field vector,
$\mathbf{H}=(H^x,H^y,H^z)$ is the magnetic field vector,
and~$\mathbf{i}=(i^x,i^y,i^z)$ is the electric current.
The conductivity~$\sigma$ is some function of~$p$ and~$\rho$, i.e., $\sigma=\sigma(p,\rho)$.

\medskip

Following~\cite{bk:SamarskyPopov_book[1992]} we firstly consider the
simplest case of one-component electric and magnetic fields. Here we
also introduce the notation and some basic concepts. In the next
section some generalizations are considered, including the case of
infinite conductivity.

\begin{figure}[ht]
\begin{adjustwidth}{-\extralength}{0cm}
\centering
\includegraphics[width=0.3\linewidth]{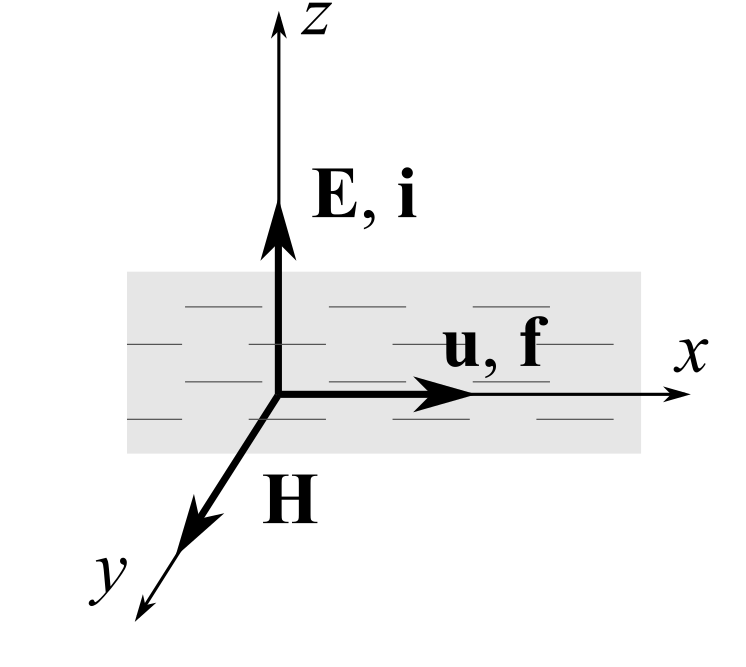}
\end{adjustwidth}
\caption{A plane one-dimensional flow for the chosen coordinates}
\label{fig:plainflow}
\end{figure}

For simplicity, the longitudinal component of the magnetic field~$\mathbf{H}$
is set to zero, and the coordinate system is chosen in such a way that $\mathbf{H} = (0, H, 0)$.
Consequently, the electric current~$\mathbf{i}$ and the electric field~$\mathbf{E}$ are also one-component vectors, i.e.,
$\mathbf{i} = (0,0, i)$, $\mathbf{E} = (0, 0, E)$. Electromagnetic force $\mathbf{f}=(f,0,0)$ acts in the $x$-direction,
and the velocity is~$\mathbf{u}=(u,0,0)$ (see Figure~\ref{fig:plainflow}).

\medskip

Given the above, the system of the one-dimensional MHD equations
with a finite conductivity $|\sigma| < \infty$ in mass Lagrangian
coordinates can be written as~\cite{bk:SamarskyPopov_book[1992]}
\begin{subequations}
\label{eq:sys0}
\begin{equation}
    \left(\frac{1}{\rho}\right)_t = u_s,
\end{equation}
\begin{equation}
\label{eq:velocity}
    u_t = -p_s + f,
    \qquad
    f = -i H / \rho,
\end{equation}
\begin{equation}
    \displaystyle
    \left(\frac{H}{\rho}\right)_t = E_s,
\end{equation}
\begin{equation}
    \displaystyle
    i = \sigma E = \kappa{\rho}H_s,
\end{equation}
\begin{equation}
\label{eq:energy}
    \varepsilon_t = -p u_s + q,
    \qquad
    q = i E / \rho,
\end{equation}
\begin{equation}
    x_t = u,
    \qquad
    x_s = 1/\rho,
\end{equation}
\end{subequations}
where $\kappa=1/(4\pi)$ and~$q$ is Joule heating per unit mass.

\medskip

In particular, we consider a polytropic gas for that the following relation holds
\begin{equation} \label{energyPolyGas}
\varepsilon = \frac{1}{\gamma - 1}\frac{p}{\rho}, \qquad \gamma = \text{const} > 1.
\end{equation}

Equation (\ref{eq:energy}) for the energy evolution can
be rewritten in the semi-divergent form
\begin{equation}
\displaystyle
\left( \varepsilon + \frac{u^2}{2} \right)_t = -(p u)_s + f u + q,
\end{equation}
or in the divergent form
\begin{equation}
\displaystyle
\left( \varepsilon + \frac{u^2}{2} + \kappa\frac{H^2}{2\rho}\right)_t
= -\left[ \left( p + \kappa \frac{H^2}{2}\right) u \right]_s + \kappa\left(E H\right)_s.
\end{equation}

Notice that the electromagnetic force $f=-i H/\rho$ can be represented in the divergent form
$f = -\kappa(H^2/2)_s$, and equation~(\ref{eq:velocity}) can be rewritten as
\begin{equation}
u_t = -\left(p + \kappa\frac{H^2}{2}\right)_s.
\end{equation}

\medskip

Further we assume $\kappa=1$ since it can be discarded by means of the scaling transformation
\begin{equation}
\tilde{s} = \kappa s,
\quad
\tilde{p} = \kappa p,
\quad
\tilde{\rho} = \kappa \rho,
\quad
\tilde{\sigma} = \kappa \sigma.
\end{equation}

\subsection{Conservative Samarskiy--Popov's schemes for system (\ref{eq:sys0})}

The family of Samarskiy--Popov's conservative difference schemes for system~(\ref{eq:sys0}) is
\begin{subequations}
\label{eq:scheme0}
\begin{equation}
\left(\frac{1}{\rho}\right)_t = u_s^{(0.5)},
\end{equation}
\begin{equation}
\displaystyle
u_t = -p^{(\alpha)}_{\bar{s}} + f,
\qquad
f = -\left( 
\frac{H\hat{H}}{2}\right)_{\bar{s}}
    = -\frac{1}{2}[\hat{i}H_*/\hat{\rho}_* + i \hat{H}_* / \rho_*],
\end{equation}
\begin{equation}
\left(\frac{H}{\rho}\right)_t = E_s^{(\beta)},
\end{equation}
\begin{equation}
i = \sigma_* E = 
    \rho_* H_{\bar{s}},
\end{equation}
\begin{equation}
\label{eq:scm0_energy}
\displaystyle
\varepsilon_t = -p^{(\alpha)} u_s^{(0.5)} + q,
\qquad
q = \frac{1}{2}[(i/\rho_*)^{(0.5)} E^{(\beta)} + (i_+ / (\rho_+)_*)^{(0.5)} E_+^{(\beta)}],
\end{equation}
\begin{equation}
x_t = u^{(0.5)},
\qquad
x_s = \frac{1}{\rho},
\end{equation}
\end{subequations}
where $0 \leqslant (\alpha, \beta) \leqslant 1$ are free parameters.
For arbitrary $\alpha$ and $\beta$ scheme (\ref{eq:scheme0}) approximates system~(\ref{eq:sys0}) up to~$O(\tau + h^2)$,
and for $\alpha = \beta = 0.5$ the scheme is of order~$O(\tau^2 + h^2)$.

Here and further~$\phi_t$, $\phi_{\check{t}}$ and~$\phi_s$,
$\phi_{\bar{s}}$ denote finite-difference derivatives of some
quantity~$\phi=\phi(t_n, s_m, u^n_m, ...)$
\begin{equation}
\def\arraystretch{2.0}
\begin{array}{c}
\displaystyle
\phi_t = \frac{\dtp{S}(\phi) - \phi}{\tau_n},
\qquad
\phi_s = \frac{\dsp{S}(\phi) - \phi}{h_m},
\\
\displaystyle
\phi_{\check{t}} = \frac{\phi - \dtm{S}(\phi)}{\tau_{n-1}},
\qquad
\phi_{\bar{s}} = \frac{\phi - \dsm{S}(\phi)}{h_{m-1}},
\end{array}
\end{equation}
which are defined with the help of the finite-difference right and
left shifts along the time and space axes correspondingly
\[
    \def\arraystretch{1.75}
    \begin{array}{c}
    \displaystyle
    \underset{\pm\tau}{S}(\phi(t_n, s_m, u^{n}_{m}, ...))
        = \phi(t_{n \pm 1}, s_m, u^{n \pm 1}_{m}, ...),
    \\
    \displaystyle
    \underset{\pm{s}}{S}(\phi(t_n, s_m, u^{n}_{m}, ...))
        = \phi(t_n, s_{m\pm 1}, u^{n}_{m \pm 1}, ...).
    \end{array}
\]
The indices~$n$ and~$m$ are respectively changed along time and space axes~$t$
and~$s$. The time and space steps are defined as follows
\begin{equation}
\def\arraystretch{1.25}
\begin{array}{c}
\displaystyle
\tau_n = \hat{\tau} = t_{n+1} - t_n = \hat{t} - t,
\qquad
\tau_{n-1} = \check{\tau} = t_n - t_{n-1} = t - \check{t},
\\
\displaystyle
h_m = h_+ = s_{m+1} - s_m = s_+ - s,
\qquad
h_{m-1} = h_- = s_{m} - s_{m-1} = s - s_-.
\end{array}
\end{equation}

\medskip

Following the Samarskiy--Popov notation throughout the text we denote
\begin{equation}
\dsp{S}(\phi) = \phi_+,
\quad
\dsm{S}(\phi) = \phi_-,
\quad
\dtp{S}(\phi) = \hat{\phi},
\quad
\dtm{S}(\phi) = \check{\phi},
\end{equation}

\begin{equation}
{\phi}^{(\alpha)}=\alpha\hat{{\phi}}+(1-\alpha){\phi}
\end{equation}
and
\begin{equation}
\label{eq:star}
{\phi}_*=({\phi}_*)^j_i=\frac{h_i {\phi}^j_{i-1/2}+h_{i-1}{\phi}^j_{i+1/2}}{h_i+h_{i-1}}.
\end{equation}
Notice that on a uniform lattice $h_i = h = \text{const}$ in its integral nodes~(\ref{eq:star}) becomes
\begin{equation}
\displaystyle
{\phi}_* = \frac{{\phi}_- + {\phi}}{2}.
\end{equation}

\begin{Remark}
The energy equation~(\ref{eq:scm0_energy}) can be reduced to one of the three following forms~\cite{bk:SamarskyPopov_book[1992]}
using equivalent algebraic transformations:
\begin{equation}
\displaystyle
\varepsilon_t = -p^{(\alpha)} u_s^{(0.5)} + q,
\end{equation}
\begin{equation}
\displaystyle
\left(\varepsilon + \frac{u^2 + u_+^2}{4} \right)_t
= -\left(p_*^{(\alpha)}u^{(0.5)}\right)_s
+ \frac{1}{2}\left[ f u^{(0.5)} + f_+ u_+^{(0.5)}\right] + q,
\end{equation}
\begin{equation}
\left(\varepsilon + \frac{u^2 + u_+^2}{4} + 
\frac{H^2}{2\rho}\right)_t
+ \left[
    \left( p_*^{(\alpha)} + 
    \frac{(H\hat{H})_*}{2}\right) u^{(0.5)}
    - E^{(\beta)} H_*^{(0.5)}
\right]_s = 0.
\end{equation}
These different forms of equation reflect the balance of certain types of energy,
i.e. they express the different physical aspects of energy conservation.
To emphasize this property, such schemes are also called \emph{completely} conservative.
\end{Remark}

\subsection{Invariance of Samarskiy--Popov's schemes}

System (\ref{eq:sys0}) can be rewritten in the following form that is more convenient for symmetry analysis
\begin{subequations}
\label{eq:sys1}
\begin{equation}
    \left(\frac{1}{\rho}\right)_t = u_s,
\end{equation}
\begin{equation}
    u_t = -\left(p + \frac{H^2}{2}\right)_s,
\end{equation}
\begin{equation}
    \left(\frac{H}{\rho}\right)_t = E_s,
\end{equation}
\begin{equation}
    \displaystyle
    \sigma E = 
    \rho H_s,
\end{equation}
\begin{equation} \label{energyEqPoly}
    p_t = -\gamma \rho p u_s + (\gamma-1) \sigma E^2,
\end{equation}
\begin{equation}
    x_t = u,
    \qquad
    x_s = 1/\rho.
\end{equation}
\end{subequations}

\begin{Remark}
Notice that for the polytropic gas with the state equation~(\ref{energyPolyGas}) one can rewrite the energy evolution equation (\ref{eq:scm0_energy}) of the Samarskiy--Popov scheme as
\begin{equation} \label{SamPopenergyPoly}
\displaystyle
p_t = - \hat{\rho} (p +  (\gamma-1)p^{(\alpha)}) u^{(0.5)}
+ (\gamma-1)\frac{\hat{\rho}}{2}
\left[
    \left(\frac{\sigma_* E}{\rho_*}\right)^{(0.5)} E^{(\beta)}
    + \left(\frac{(i_+)_* E_+ }{ (\rho_+)_*}\right)^{(0.5)} E_+^{(\beta)}
\right].
\end{equation}
In this form the energy evolution equation corresponds to equation~(\ref{energyEqPoly}).
\end{Remark}

\medskip

Calculations show~\cite{bk:DorKozMelKap_PlainFlows_2021} that the Lie algebra admitted by the system for an arbitrary $\sigma=\sigma(p,\rho)$ is\footnote{Here and further the notation~$\displaystyle \partial_f \equiv \frac{\partial}{\partial{f}}$ is used.}
\begin{equation} \label{kern0}
\def\arraystretch{1.75}
\begin{array}{c}
X_1 = \partial_{t},
\qquad
X_2 = \partial_{s},
\qquad
X_3 = \partial_{x},
\qquad
X_4 = t \partial_{x} + \partial_{u}.
\end{array}
\end{equation}

The group generator
\begin{equation} \label{Xgen}
X = \xi^t \partial_t
+ \xi^s \partial_s
+ \eta \partial_x
\end{equation}
is prolonged to the finite-difference space as follows~\cite{Dor_1,bk:Dorodnitsyn[2011]}
\begin{equation}
\displaystyle
\tilde{X} = \sum_{k,l=-\infty}^{\infty} \dtm{S^k}\dsm{S^l}(X),
\end{equation}
The scheme of the form
\begin{subequations} \label{ScmGenAll}
\begin{equation}\label{ScmGen}
\Phi(t,\check{t}, \hat{t}, s, s_+, s_-, u, u_+, u_-, \hat{u}, \check{u}, \hat{u}_+, \check{u}_+, \hat{u}_-, \check{u}_-, ...) = 0,
\end{equation}
\begin{equation}\label{ScmGenMesh}
\displaystyle h_+ = h_- = h, \qquad
\hat{\tau} = \check{\tau} = \tau, \qquad  (\vec{\tau}, \vec{h})=0
\end{equation}
\end{subequations}
defined on a uniform orthogonal mesh is \emph{invariant}
if the following criterion of invariance holds~\cite{bk:Dorodnitsyn[2011]}
\begin{equation} \label{SchInvCond}
\def\arraystretch{1.5}
\begin{array}{l}
\displaystyle
\tilde{X}\Phi|_{(\ref{ScmGenAll})} = 0,
\\
\displaystyle
\tilde{X}(\hat{\tau} - \check{\tau})|_{(\ref{ScmGenAll})} = 0,
\quad
\tilde{X}(h_+ - h_-)|_{(\ref{ScmGenAll})} = 0.
\end{array}
\end{equation}
To preserve uniformness and orthogonality of the mesh it is also needed~\cite{Dor_1,bk:Dorodnitsyn[2011]}
\begin{equation} \label{mesh_conds_uni}
  \dsp{D}\dsm{D}(\xi^s) = 0,
  \qquad
  \dtp{D}\dtm{D}(\xi^t) = 0,
\end{equation}
\begin{equation} \label{mesh_conds_ortho}
  \underset{\pm{s}}{D}(\xi^t) = -\underset{\pm{\tau}}{D}(\xi^s),
\end{equation}
where $\underset{\pm\tau}{D}$ and $\underset{\pm{s}}{D}$ are
finite-difference differentiation operators
\[
    \underset{+\tau}{D} = \frac{\underset{+\tau}{S} - 1}{\tau_n},
    \quad
    \underset{-\tau}{D} = \frac{1 - \underset{-\tau}{S}}{\tau_{n-1}},
    \quad
    \underset{+s}{D} = \frac{\underset{+s}{S} - 1}{h_m},
    \quad
    \underset{-s}{D} = \frac{1 - \underset{-s}{S}}{h_{m-1}}.
\]

\medskip

One can verify that scheme~(\ref{eq:scheme0}) is indeed invariant
with respect to the generators $X_1$, ..., $X_4$ and all the
generators~(\ref{kern0}) satisfy the mesh orthogonality and
uniformness
conditions~(\ref{mesh_conds_uni}),~(\ref{mesh_conds_ortho}). Hence,
one can use an orthogonal uniform mesh~(\ref{ScmGenMesh}) which is
invariant one.

\subsection{Conservation laws possessed by Samarskiy--Popov's scheme}

All the conservation laws of system~(\ref{eq:sys1})
have their finite-difference counterparts for the Samarskiy--Popov scheme.
They are given in Table~\ref{tab:CLs-table}.

\medskip

An additional conservation law
\begin{equation} \label{CLadd}
\left(\frac{s H}{\rho}\right)_t + (H - s E)_s = 0.
\end{equation}
only occurs in case~$\sigma = \rho$~\cite{bk:DorKozMelKap_PlainFlows_2021}.
In this case, system~(\ref{eq:sys1}) admits two more symmetries, namely
\begin{equation}
\def\arraystretch{1.5}
\begin{array}{c}
X_5 = s \partial_{s}
    - x \partial_{x}
    - u \partial_{u}
    + 2\rho \partial_{\rho}
    - E \partial_{E},
\\
X_6 = 2 t\partial_{t} + 2 s \partial_{s}
    - 2 u \partial_{u}
    + 2 p \partial_{p}
    + 2\rho \partial_{\rho}
    - 3 E \partial_{E}
    - H \partial_{H}.
\end{array}
\end{equation}
Conservation law~(\ref{CLadd}) has its finite-difference counterpart which can be found by direct calculations.

\medskip

Notice that neither conservation law~(\ref{CLadd}) nor the center-of-mass law
\begin{equation} \label{CLcmass}
\left(t u - x\right)_t + \left(t \left(p + \frac{H^2}{2}\right)\right)_s = 0
\end{equation}
was mentioned in~\cite{bk:SamarskyPopov_book[1992]}.
Perhaps, the authors of~\cite{bk:SamarskyPopov_book[1992]} have known the finite-difference
analogue of~(\ref{CLcmass}).

\begingroup
\def\arraystretch{2.0}
\begin{table}[H]
\caption{Differential and difference conservation laws\label{tab:CLs-table}}
\begin{adjustwidth}{-\extralength}{0cm}
\newcolumntype{C}{>{\centering\arraybackslash}X}
\begin{tabularx}{\fulllength}{cCCC}
\toprule
\textbf{\#}
&
\textbf{Conservation laws of system}~(\protect\ref{eq:sys1})
&
\textbf{Conservation laws of scheme}~(\protect\ref{eq:scheme0})
&
\textbf{Physics interpretation}\\
\midrule
\multicolumn{4}{c}{$\sigma= \sigma(p, \rho)$}\\
\midrule
1
&
$
\left(\frac{1}{\rho}\right)_t - u_s = 0
$
&
$
\left(\frac{1}{\rho}\right)_t  - \left(u^{(0.5)}\right)_s = 0
$
&
Mass conservation
\\ 
2
&
$
\left(\frac{H}{\rho}\right)_t - E_s = 0
$
&
$
\left(\frac{H}{\rho}\right)_t - (E^{(\beta)})_s = 0
$
& Magnetic flux conservation
\\ 
3
&
$
u_t + \left(p + \frac{H^2}{2}\right)_s = 0
$
&
$
u_t + \left(p_-^{(\alpha)} + 
\frac{H_-\hat{H}_-}{2} \right)_s = 0
$
& Momentum conservation
\\ 
4
&
$
\left(t u - x\right)_t + \left(t \left(p + \frac{H^2}{2}\right)\right)_s = 0
$
&
$
\left({\displaystyle \frac{t + \check{t}}{2}} \, u - x\right)_t +~\left[ t\left(p_-^{(\alpha)} + 
\frac{H_-\hat{H}_-}{2} \right)\right]_s~=~0
$
& Center of mass law
\\ 
5
&
$
\begin{array}{l}
\left(
    \varepsilon
    + \frac{u^2}{2}
    + \frac{H^2}{2\rho}
\right)_t + \\
\quad + \left(
    \left(p + \frac{H^2}{2}\right) u
    -E H
\right)_s = 0
\end{array}
$
&
$
\begin{array}{l}
\left(\varepsilon + \frac{u^2 + u_+^2}{4} + 
\frac{H^2}{2\rho}\right)_t + \\
\quad + \left[
    \left( p_*^{(\alpha)} + 
    \frac{(H\hat{H})_*}{2}\right) u^{(0.5)}
    - 
    E^{(\beta)} H_*^{(0.5)}
\right]_s = 0
\end{array}
$
&
Energy conservation
\\
 \midrule
\multicolumn{4}{c}{$\sigma= \rho$}
\\
 \midrule
6
&
$
\left(\frac{s H}{\rho}\right)_t + (H - s E)_s = 0
$
&
$
\left(\frac{s H}{\rho}\right)_t + (H_-^{(\beta)} - s_- E^{(\beta)})_s = 0
$
&
\unknownCL
\\
\bottomrule
\end{tabularx}
\end{adjustwidth}
\end{table}
\endgroup

\section{Generalizations of the Samarskiy--Popov schemes for MHD equations}
\label{sec:general}

\subsection{The case of finite conductivity}
\label{sec:genFiniteCond}

We consider a more general case
$\mathbf{H} = (H^0, H^y, H^z)$,
$\mathbf{E} = (0, E^y, E^z)$,
$\mathbf{i} = (0, i^y, i^z)$,
$\mathbf{u}=(u,v,w)$, $\mathbf{x}=(x,y,z)$, and $H^0 = \text{const}$.
Here we used the fact that the coordinate system can always be chosen
in such a way that the first component of the vector field~$\mathbf{E}$ is equal to zero.

\medskip

Further we consider equations~(\ref{eq:sys1ext}), where,
by analogy with~(\ref{eq:sys1}), the energy evolution equation~(\ref{eq:energyVec}) is written as
\begin{equation}
p_t = -\gamma p \rho u_{{s}} + (\gamma - 1) \sigma ((E^y)^2 + (E^z)^2),
\end{equation}

\medskip

A generalization of scheme~(\ref{eq:scheme0}) for $\mathbf{E}=(0, E^y, 0)$ and $\mathbf{H}=(H^0, 0, H^z)$ is given in~\cite{bk:SamarskyPopov_book[1992]}.
Since the MHD equations are almost symmetric in terms of the components $E^y$, $E^z$ and $H^y$, $H^z$, one can extend the scheme proposed in~\cite{bk:SamarskyPopov_book[1992]} as follows
\begin{subequations}
\label{eq:scheme0ext}
\begin{equation}
\left(\frac{1}{\rho}\right)_t = u_s^{(0.5)},
\end{equation}
\begin{equation}
\displaystyle
u_t = -\left( p^{(\alpha)}
    + \frac{H^y\hat{H}^y + H^z\hat{H}^z} {2}\right)_{\bar{s}},
\qquad
v_t = H^0 (H^y)^{(0.5)}_{\bar{s}},
\qquad
w_t = H^0 (H^z)^{(0.5)}_{\bar{s}},
\end{equation}
\begin{equation}
\left(\frac{H^y}{\rho}\right)_t = H^0 v^{(0.5)}_s + (E^z)_s^{(\beta_1)},
\qquad
\left(\frac{H^z}{\rho}\right)_t = H^0 w^{(0.5)}_s - (E^y)_s^{(\beta_2)},
\end{equation}
\begin{equation}
i^y = \sigma_* E^y = -\rho_* {H}^z_{\bar{s}},
\qquad
i^z = \sigma_* E^z = \rho_* {H}^y_{\bar{s}},
\end{equation}
\begin{equation}
\label{eq:scm1_energy}
\displaystyle
\varepsilon_t = -p^{(\alpha)} u_s^{(0.5)} + q^y + q^z,
\end{equation}
\begin{equation}
x_t = u^{(0.5)},
\qquad
y_t = v^{(0.5)},
\qquad
z_t = w^{(0.5)},
\qquad
x_s = \frac{1}{\rho},
\end{equation}
\end{subequations}
where
$0 \leqslant (\alpha, \beta_1, \beta_2) \leqslant 1$
and
\[
q^y = \frac{1}{2}[(i^y/\rho_*)^{(0.5)} (E^y)^{(\beta_2)} + (i^y_+ / (\rho_+)_*)^{(0.5)} (E^y)_+^{(\beta_2)}],
\]
\[
q^z = \frac{1}{2}[(i^z/\rho_*)^{(0.5)} (E^z)^{(\beta_1)} + (i^z_+ / (\rho_+)_*)^{(0.5)} (E^z)_+^{(\beta_1)}].
\]

Notice that this generalization of the scheme was discussed in~\cite{bk:SamarskyPopov_book[1992]} but it was not given explicitly.

\medskip

\begin{Remark}
One can generalize~(\ref{SamPopenergyPoly}) for scheme~(\ref{eq:scheme0ext}), (\ref{energyPolyGas}) as follows
\begin{equation}
p_t + \hat{\rho} u_s^{(0.5)} (p + (\gamma-1) p^{(\alpha)}) + (1-\gamma) \hat{\rho}(q^y + q^z) = 0.
\end{equation}
\end{Remark}

\medskip

According to~\cite{bk:DorKozMelKap_PlainFlows_2021}, the symmetries admitted by system~(\ref{eq:sys1ext})
are the following.

\begin{enumerate}

\item
If $H^0 = 0$ and $\sigma$ is arbitrary, then the admitted Lie algebra is
\begin{equation}
\def\arraystretch{1.75}
\begin{array}{c}
{X_1} = \partial_{t},
\qquad
{X_2} = \partial_{s},
\qquad
{X_3} = \partial_{x},
\qquad
{X_4} = t \partial_{x} + \partial_{u},
\\
{X_5} = E^z \partial_{E^y} - E^y \partial_{E^z} + H^z \partial_{H^y} - H^y \partial_{H^z}.
\end{array}
\end{equation}
In case $\sigma=\rho$, there are two additional generators are admitted, namely
\begin{equation}
\def\arraystretch{1.75}
\begin{array}{c}
X_{1a} = s \partial_{s}
    - x \partial_{x}
    - u \partial_{u}
    + 2\rho \partial_{\rho}
    - E^y \partial_{E^y}
    - E^z \partial_{E^z},
\\
X_{2a} =
    2 t \partial_{t}
    + 2 x \partial_{x}
    - 2 p\partial_{p}
    - 2 \rho \partial_{\rho}
    - E^y \partial_{E^y}
    - E^z \partial_{E^z}
    - H^y \partial_{H^y}
    - H^z \partial_{H^z}.
\end{array}
\end{equation}

There are also two more conservation laws in the latter case (see Table~\ref{tab:CLsExt-table}).
Additional conservation laws do not occur for any other forms of the function~$\sigma$.

\item
If $H^0 \neq 0$ and $\sigma$ is arbitrary then the admitted Lie algebra is
\begin{equation}
\def\arraystretch{1.75}
\begin{array}{c}
{X_1} = \partial_{t},
\qquad
{X_2} = \partial_{s},
\qquad
{X_3} = \partial_{x},
\qquad
{X_4} = t \partial_{x} + \partial_{u},
\\
{X_5} = E^z \partial_{E^y} - E^y \partial_{E^z} + H^z \partial_{H^y} - H^y \partial_{H^z}
    + w\partial_{v} - v\partial_{w}
    + z \partial_{y} - y \partial_{z},
\\
{X_6} = t \partial_{y} + \partial_{v},
\qquad
{X}_{7} = t \partial_{z} + \partial_{w},
\qquad
{X}_{8}= h_1(s) \partial_{y},
\qquad
{X}_{9}= h_2(s) \partial_{z},
\end{array}
\end{equation}
where $h_1$ and $h_2$ are arbitrary functions of $s$.

Additional conservation laws do not occur for any specific~$\sigma$.

\end{enumerate}

\medskip

In both the cases above, scheme~(\ref{eq:scheme0ext}) is \emph{invariant}.
The rotation generator~$X_5$ is only admitted for~$\beta_1 = \beta_2$.
The remaining generators are admitted by the scheme for any set of parameters $\alpha$, $\beta_1$, $\beta_2$ and $\gamma > 1$.

\medskip

The conservation laws possessed by system~(\ref{eq:scheme0ext}) and their finite-difference
counterparts are given in~Table~\ref{tab:CLsExt-table}.
Here and further, conservation laws whose fluxes vanish for~$H^0=0$ are marked with~$\dag$.
In case~$H^0=0$, their densities preserve along the pathlines.

\subsection{The case of infinite conductivity $\sigma \to \infty$}
\label{sec:genInfiniteCond}

In this case, system~(\ref{eq:sys1ext}) reduces to
\begin{subequations} \label{eq:sys1ext_inf}
\begin{equation} \label{MassEq}
\left(\frac{1}{\rho}\right)_t = u_s,
\end{equation}
\begin{equation} \label{VelosityEq}
u_t = -\left( p + \frac{(H^y)^2 + (H^z)^2}{2} \right)_s,
\qquad
x_t = u,
\end{equation}
\begin{equation}
v_t = H^0 H^y_{{s}},
\qquad
y_t = v,
\end{equation}
\begin{equation}
w_t = H^0 H^z_{{s}},
\qquad
z_t = w,
\end{equation}
\begin{equation}
\left(\frac{H^y}{\rho}\right)_t = (H^0 v)_s,
\end{equation}
\begin{equation}
\left(\frac{H^z}{\rho}\right)_t = (H^0 w)_s,
\end{equation}
\begin{equation} \label{polyEnergyEq}
p_t = -\gamma p \rho u_{s},
\end{equation}
\end{subequations}
where the internal energy is given by~(\ref{energyPolyGas}).

In addition to the analogues of conservation laws presented in the previous section, system~(\ref{eq:sys1ext_inf})
possesses the conservation law of angular momentum, namely
\begin{equation}
\label{clMom}
(z v - y w)_t + \left( H^0 (y H^z - z H^y) \right)_s = 0.
\end{equation}

\bigskip

As $\sigma \to \infty$, scheme (\ref{eq:scheme0ext}) becomes
\begin{subequations} \label{SamPopScmInfMod0}
\begin{equation}
\left(\frac{1}{\rho}\right)_t = u_s^{(0.5)},
\end{equation}
\begin{equation}
\displaystyle
u_t = -\left( p^{(\alpha)}
    + \frac{H^y\hat{H}^y + H^z\hat{H}^z} {2}\right)_{\bar{s}},
\qquad
v_t = H^0 (H^y)^{(0.5)}_{\bar{s}},
\qquad
w_t = H^0 (H^z)^{(0.5)}_{\bar{s}},
\end{equation}
\begin{equation}
\left(\frac{H^y}{\rho}\right)_t = H^0 v^{(0.5)}_s,
\qquad
\left(\frac{H^z}{\rho}\right)_t = H^0 w^{(0.5)}_s,
\end{equation}
\begin{equation} \label{schemeinfEnergyPressure}
\displaystyle
\varepsilon_t = -p^{(\alpha)} u_s^{(0.5)},
\end{equation}
\begin{equation}
x_t = u^{(0.5)},
\qquad
y_t = v^{(0.5)},
\qquad
z_t = w^{(0.5)},
\qquad
x_s = \frac{1}{\rho}.
\end{equation}
\end{subequations}

\begingroup
\def\arraystretch{2.0}
\begin{table}[H]
\caption{Differential and difference conservation laws for the extended scheme\label{tab:CLsExt-table}}
\begin{adjustwidth}{-\extralength}{0cm}
\newcolumntype{C}{>{\centering\arraybackslash}X}
\begin{tabularx}{\fulllength}{cCCC}
\toprule
\textbf{\#}
&
\textbf{Conservation laws of system}~(\protect\ref{eq:sys1ext})
&
\textbf{Conservation laws of scheme}~(\protect\ref{eq:scheme0ext})
&
\textbf{Physics interpretation}
\\
\midrule
\multicolumn{4}{c}{$\sigma= \sigma(p, \rho)$}
\\
\midrule
1
&
$
\left(\frac{1}{\rho}\right)_t - u_s = 0
$
&
$
\left(\frac{1}{\rho}\right)_t  - \left(u^{(0.5)}\right)_s = 0
$
&
Mass conservation
\\ 
2
&
$
\left(\frac{H^y}{\rho}\right)_t - (E^z + H^0 v)_s = 0
$
&
$
\left(\frac{H^y}{\rho}\right)_t - \left((E^z)^{(\beta_1)} + H^0 v^{(0.5)}\right)_s = 0
$
&
Magnetic flux conservation
\\ 
3
&
$
\left(\frac{H^z}{\rho}\right)_t + (E^y - H^0 w)_s = 0
$
&
$
\left(\frac{H^z}{\rho}\right)_t + \left((E^y)^{(\beta_2)} - H^0 w^{(0.5)}\right)_s = 0
$
&
Magnetic flux conservation
\\ 
4
&
$
u_t + \left(p + \frac{(H^y)^2+ (H^z)^2}{2}\right)_s = 0
$
&
$
u_t + \left(p_-^{(\alpha)} + 
\frac{H^y_-\hat{H}^y_- + H^z_-\hat{H}^z_-}{2} \right)_s = 0
$
&
Momentum conservation
\\ 
$5^\dag$
&
$
v_t - (H^0 H^y)_s = 0
$
&
$
v_t - \left(H^0 (H^y_-)^{(0.5)}\right)_{s} = 0
$
&
Momentum conservation
\\ 
$6^\dag$
&
$
w_t - (H^0 H^z)_s = 0
$
&
$
w_t - \left(H^0 (H^z_-)^{(0.5)}\right)_{s} = 0
$
&
Momentum conservation
\\ 
7
&
$
\left(t u - x\right)_t + \left[t \left(p + \frac{(H^y)^2+(H^z)^2}{2}\right)\right]_s = 0
$
&
$
\left({\displaystyle \frac{t + \check{t}}{2}} \, u - x\right)_t
+ \left[ t\left(p_-^{(\alpha)} + 
\frac{H^y_-\hat{H}^y_- + H^z_-\hat{H}^z_-}{2} \right)\right]_s = 0
$
&
Center of mass law
\\ 
$8^\dag$
&
$
(t v - y)_t - (t H^0 H^y)_s = 0
$
&
$
\left({\displaystyle \frac{t + \check{t}}{2}} \, v - y\right)_t
- \left( t H^0 (H^y_-)^{(0.5)} \right)_s = 0
$
&
Center of mass law
\\ 
$9^\dag$
&
$
(t w - z)_t - (t H^0 H^z)_s = 0
$
&
$
\left({\displaystyle \frac{t + \check{t}}{2}} \, w - z\right)_t
- \left( t H^0 (H^z_-)^{(0.5)} \right)_s = 0
$
&
Center of mass law
\\ 
10
&
$
{\def\arraystretch{2}
\begin{array}{l}
\left(
    \varepsilon
    + \frac{u^2 + v^2 + w^2}{2}
    + \frac{(H^y)^2 + (H^z)^2}{2\rho}
\right)_t +
\\
+ \left[
    \left(p + \frac{(H^y)^2 + (H^z)^2}{2}\right) u
    + E^y H^z - E^z H^y
    \right.
\\
    \left.
    -H^0 (v H^y + w H^z)
\right]_s = 0
\end{array}
} 
$
&
$
{\def\arraystretch{2}
\begin{array}{l}
\left(\varepsilon + \frac{u^2 + u_+^2 + v^2 + v_+^2 + w^2 + w_+^2}{4} + 
\frac{(H^y)^2 + (H^z)^2}{2\rho}\right)_t + \\
+ \left[
    \left( p_*^{(\alpha)} + 
    \frac{(H^y\hat{H}^y + H^z\hat{H}^z)_*}{2}\right) u^{(0.5)}
    \right.
    \\
    \left.
    + (E^y)^{(\beta_2)} (H^z_*)^{(0.5)}
    - (E^z)^{(\beta_1)} (H^y_*)^{(0.5)}
    \right.
    \\
    \left.
    - H^0 (
        v^{(0.5)}(H^y_*)^{(0.5)}
        + w^{(0.5)}(H^z_*)^{(0.5)}
    )
\right]_s = 0
\end{array}
} 
$
&
Energy conservation
\\
\midrule
\multicolumn{4}{c}{$H^0 = 0$}
\\
\midrule
$11^\dag$
&
$
\left(z v - y w\right)_t = 0
$
&
$
\displaystyle
({z} {v}^{(0.5)} - {y} {w}^{(0.5)})_{\check{t}} = 0
$
&
Angular momentum conservation
\\
\midrule
\multicolumn{4}{c}{$H^0 = 0$, $\sigma= \rho$}
\\
\midrule
12
&
$
\left(\frac{s H^y}{\rho}\right)_t + (H^y - s E^z)_s = 0
$
&
$
\left(\frac{s H^y}{\rho}\right)_t + ((H_-^y)^{(\beta_1)} - s_- (E^z)^{(\beta_1)})_s = 0
$
&
\unknownCL
\\
13
&
$
\left(\frac{s H^z}{\rho}\right)_t + (H^z + s E^y)_s = 0
$
&
$
\left(\frac{s H^z}{\rho}\right)_t + ((H_-^z)^{(\beta_2)} + s_- (E^y)^{(\beta_2)})_s = 0
$
&
\unknownCL
\\
\bottomrule
\end{tabularx}
\end{adjustwidth}
\end{table}
\endgroup


One can verify that scheme~(\ref{SamPopScmInfMod0}) is invariant one.
As the symmetries of~(\ref{eq:sys1ext_inf}) and the corresponding difference schemes
are reviewed in Section~\ref{sec:isentrop}, we defer our discussion until then.

\subsubsection{Conservation of angular momentum and energy}

Apparently, the latter scheme does not preserve angular momentum, i.e. it does not possess a difference analogue of the conservation law~(\ref{clMom}). One can verify it by algebraic manipulations with the scheme
or with the help of the finite-difference analogue of the direct method~\cite{bk:ChevDorKap2020,dorodnitsyn2019shallow}.
We overcome this issue by modifying the latter scheme as follows
\begin{subequations} \label{SamPopScmInfMod1}
\begin{equation}
\left(\frac{1}{\rho}\right)_t = u_s^{(0.5)},
\end{equation}
\begin{equation}
\displaystyle
u_t = -\left( p^{(\alpha)}
    + \frac{H^y\hat{H}^y + H^z\hat{H}^z} {2}\right)_{\bar{s}},
\qquad
v_t = H^0 \hat{H}^y_{\bar{s}},
\qquad
w_t = H^0 \hat{H}^z_{\bar{s}},
\end{equation}
\begin{equation} \label{usws}
\left(\frac{H^y}{\rho}\right)_t = H^0 v_s,
\qquad
\left(\frac{H^z}{\rho}\right)_t = H^0 w_s,
\end{equation}
\begin{equation}
\displaystyle
\varepsilon_t = -p^{(\alpha)} u_s^{(0.5)},
\end{equation}
\begin{equation}
x_t = u^{(0.5)},
\qquad
y_t = v,
\qquad
z_t = w,
\qquad
x_s = \frac{1}{\rho}.
\end{equation}
\end{subequations}
This allows one to obtain the whole set of finite-difference analogues of the conservation laws of equation~(\ref{eq:sys1ext_inf}) excluding the conservation of the entropy along the pathlines.
The conservation laws are presented in Table~\ref{tab:CLsInf}.
Notice that the three-layer conservation law of energy given in the table can be rewritten in the following two-layer form by means of~(\ref{usws})
\begin{multline}
\left[\varepsilon + \frac{u^2 + u_+^2 + v^2 + v_+^2 + w^2 + w_+^2}{4}
+ \frac{(H^y)^2 + (H^z)^2}{2{\rho}}
+\frac{\tau}{2} H^0 (H^y v_s + H^z w_s)
\right]_t
\\
+ \left[
    \left( {p}_*^{(\alpha)} +
    \frac{(H^y\hat{H}^y + H^z\hat{H}^z)_*}{2} \right) {u}^{(0.5)}
    - H^0 (
        {v}^{(0.5)} \hat{H}^y_*
        + {w}^{(0.5)} \hat{H}^z_*
    )
\right]_s = 0
\end{multline}
Also, in order to verify the conservation law~(\ref{clMom}), one has to consider the following equations
which can be obtained by integration of~(\ref{usws})
\begin{equation}
y_s = \frac{H^y}{H^0 \rho},
\qquad
z_s = \frac{H^z}{H^0 \rho}.
\end{equation}

We also notice that the modified scheme~(\ref{SamPopScmInfMod1}) is still invariant and a completely conservative one.

\subsubsection{Conservation of the entropy along the pathlines}

From the latter system~(\ref{eq:sys1ext_inf}) it follows
\begin{equation}
\left(\frac{p}{\rho^\gamma}\right)_t = S_t = 0.
\end{equation}
This represents the conservation of the entropy~$S$ along pathlines
which is a crucial difference between the finite and infinite
conductivity cases.

\medskip

It is known~\cite{DORODNITSYN2019201} that the Samarskiy--Popov scheme for polytropic gas
does not preserve the entropy~$S$ for arbitrary~$\gamma$.
However, the following relation holds on solutions of the system
\begin{equation} \label{entropyDfferential}
\frac{\hat{p} - p}{p^{(\alpha)}} = \gamma \frac{\hat{\rho} - \rho}{\rho^{(\alpha)}},
\end{equation}
which approximates the differential relation
\begin{equation}
\frac{dp}{p} = \gamma \frac{d\rho}{\rho}.
\end{equation}
The latter relation holds along trajectories of the particles up to~$O(\tau)$ for~$\alpha\neq 0.5$
or up to~$O(\tau^2)$ for~$\alpha = 0.5$.

In~\cite{DORODNITSYN2019201}, an entropy preserving invariant scheme for gas dynamics equations in case of polytropic gas with $\gamma=3$ was proposed. This scheme conserves the entropy along the pathlines but has only one conservation law, namely the conservation law of mass. It seems that the conservation of entropy by the difference scheme usually leads
to the ``loss'' of some other conservation laws.

\bigskip

Here we propose a way of preserving the entropy along the pathlines for polytropic gas with integer values of adiabatic exponent~$\gamma \geqslant 2$ for scheme~(\ref{SamPopScmInfMod1}). We show that this can be done by choosing appropriate approximations of the state equation~(\ref{energyPolyGas}).

\smallskip

We notice that by means of (\ref{energyPolyGas}) and (\ref{MassEq}),
equation~(\ref{polyEnergyEq}) can be represented as the identity
\begin{equation} \label{polyGasIdentity}
\left( \frac{1}{\gamma-1} \, \rho^{\gamma-1}\right)_t = \rho^{\gamma-2} \rho_t.
\end{equation}
In the finite-difference case, the rules of differentiation are different.
As a result, not every approximation of the latter identity is a finite-difference identity.
For a proper discrete analogue of~(\ref{polyGasIdentity}) the right hand side of
the identity should also be expressed in the divergent form as well as the left hand side.
Choosing the difference approximation for the scheme in the case of a polytropic gas,
one has an additional ``degree of freedom'':
the choice of approximation for the state equation~(\ref{energyPolyGas}).
This should be done so that both the left and right hand sides of the resulting approximation for~(\ref{polyGasIdentity}) are divergent expressions.
Notice that this does not affect the conservativeness of the total energy conservation law equation
since it does not depend on any specific form of the equation of state.

\medskip

Further we consider the shifted version of the equation~(\ref{schemeinfEnergyPressure})
\begin{equation} \label{EnergyEqShifted}
\displaystyle
\check{\varepsilon}_t
= -\check{p}^{(\alpha)} \check{u}_s^{(0.5)}
= -\check{p}^{(\alpha)}  \left(\frac{1}{\check{\rho}}\right)_t.
\end{equation}
First, we choose the following approximation of the state equation~(\ref{energyPolyGas}) for $\gamma=2$,
\begin{equation} \label{stateEqapproxG=2}
\varepsilon = \frac{p^{(\alpha)}}{\hat{\rho}}.
\end{equation}
Substituting (\ref{stateEqapproxG=2}) into~(\ref{EnergyEqShifted}), one derives
\begin{equation}
\frac{p^{(\alpha)}}{\hat{\rho}} - \frac{\check{p}^{(\alpha)}}{\rho}
= -\tau \check{p}^{(\alpha)} \left(\frac{1}{\check{\rho}}\right)_t.
\end{equation}
Solving with resect to $\check{p}^{(\alpha)}$, one gets
\begin{equation}
\check{p}^{(\alpha)} = \frac{p^{(\alpha)}\check{\rho}}{\hat{\rho}}.
\end{equation}
The latter equation can be rewritten as
\begin{equation}\label{SCanBeIntegrated}
\frac{p^{(\alpha)}}{\check{p}^{(\alpha)}} = \frac{\rho\hat{\rho}}{\check{\rho}\rho}.
\end{equation}
Equation (\ref{SCanBeIntegrated}) can be integrated, i.e.,
\begin{equation}\label{schemeEntropyCLinf2}
\left(\frac{\check{p}^{(\alpha)}}{\rho\check{\rho}}\right)_t = S_t = 0.
\end{equation}
This means conservation of entropy $S$ along pathlines for $\gamma = 2$ on \emph{two} time layers.
We have achieved the integrability of the difference analogue of equation~(\ref{polyEnergyEq}) by choosing a suitable approximation for the state equation.


\medskip

In a similar way one can arrive at the conservation of entropy for~$\gamma=3$, namely
\begin{equation}
\label{schemeEntropyCLinf3}
\varepsilon = \frac{\rho p^{(\alpha)}}{\hat{\rho}(\hat{\rho} + \rho)},
\qquad
\left(\frac{2 \check{p}^{(\alpha)}}{\rho\check{\rho}(\rho + \check{\rho})}\right)_t = 0.
\end{equation}
Similarly, for $\gamma=4$
\begin{equation}
\varepsilon = \frac{\rho^2 p^{(\alpha)}}{\hat{\rho}(\hat{\rho}^2 + \rho\hat{\rho} + \rho^2)},
\qquad
\left(\frac{3 \check{p}^{(\alpha)}}{\rho\check{\rho}(\rho^2 + \rho\check{\rho}+ \check{\rho}^2)}\right)_t = 0,
\end{equation}
etc.

Thus, by induction, one establishes the following general formula for an arbitrary natural~$\gamma \geqslant 2$
\begin{equation} \label{schemeEntropyCLinfN}
\displaystyle
\varepsilon = \frac{p^{(\alpha)}}{\sum_{k=0}^{\gamma-2} \hat{\rho}^{\gamma-k-1} {\rho}^{k - \gamma + 2}},
\qquad
\left(\frac{(\gamma-1)\check{p}^{(\alpha)}}{\sum_{k=0}^{\gamma-2} \rho^{\gamma-k-1} \check{\rho}^{k+1} }\right)_t = 0.
\end{equation}
Entropy preservation formula~(\ref{schemeEntropyCLinfN}) are presented in Table~\ref{tab:CLsInf} among the other conservation laws.

\begin{Remark}
From the preservation of entropy~$S_t = 0$
in the differential case it follows (for simplicity, we consider the specific case~$\gamma = 2$)
\begin{equation} \label{entropyIntegral}
\displaystyle
\int_0^T \left(\frac{p}{\rho^2}\right)_t dt
= \frac{p(T,s)}{\rho^2(T,s)} - \frac{p(0, s)}{\rho^2(0, s)} = \text{const}.
\end{equation}
Since the constant can be omitted, this means
\[
\frac{p(0, s)}{\rho^2(0, s)} = \frac{p(T,s)}{\rho^2(T,s)}.
\]
In the finite-difference case, by means of~(\ref{schemeEntropyCLinf2}), one derives
the following analogue of~(\ref{entropyIntegral})
\[
\displaystyle
\sum_{k=0}^{N}
\left(\frac{\alpha p_m^{n+k} +(1 - \alpha)p_m^{n-1+k}}{\rho_m^{n+k}\rho_m^{n-1+k}}\right)_t \tau
= \frac{\alpha p_m^{n+N} +(1 - \alpha)p_m^{n-1+N}}{\rho_m^{n+N}\rho_m^{n-1+N}}
-\frac{\alpha p_m^{n} +(1 - \alpha)p_m^{n-1}}{\rho_m^{n}\rho_m^{n-1}} = 0,
\]
where $N = \lceil{T/\tau}\rceil$ and we recall that~$p^{(\alpha)} = \alpha \hat{p} + (1-\alpha)p$.
Similar to the differential case, the latter gives
\[
\frac{\alpha p_m^{n} +(1 - \alpha)p_m^{n-1}}{\rho_m^{n}\rho_m^{n-1}}
=
\frac{\alpha p_m^{n+N} +(1 - \alpha)p_m^{n-1+N}}{\rho_m^{n+N}\rho_m^{n-1+N}}
\]
which means entropy preservation for a given liquid particle.
\end{Remark}

\begin{Remark}
The approach described above can also lead to entropy conservation
for \emph{rational} values of~$\gamma$. Without proof  of the
existence of a general formula we present the result for~$\gamma =
5/3$ which  occurs for one-atomic ideal gas. One can verify that
for~$\gamma = 5/3$ the approximation
\begin{equation}
\varepsilon = {p}^{(\alpha)}\, \frac{\hat{\rho}^{2/3} + (\rho\hat{\rho})^{1/3} + \rho^{2/3}}{\rho^{1/3}\hat{\rho}(\hat{\rho}^{1/3} + \rho^{1/3})}
= \frac{3}{2} \frac{p}{\rho} + O(\tau)
\end{equation}
for the internal energy $\varepsilon$ leads to the following preservation of entropy
\begin{equation}
\left(
\frac{2}{3}\,
\check{p}^{(\alpha)}\, \frac{\check{\rho}^{2/3} + (\rho\check{\rho})^{1/3} + \rho^{2/3}}{\rho\check{\rho}(\check{\rho}^{1/3} + \rho^{1/3})}
\right)_t = 0.
\end{equation}
\end{Remark}

\begingroup
\def\arraystretch{2.2}
\begin{table}[H]
\caption{Differential and difference conservation laws for the modified scheme~(\ref{SamPopScmInfMod1})
for an arbitrary entropy~$S(s)$ in case $\sigma \to \infty$\label{tab:CLsInf}}
\begin{adjustwidth}{-\extralength}{0cm}
\newcolumntype{C}{>{\centering\arraybackslash}X}
\begin{tabularx}{\fulllength}{cCCC}
\toprule
\textbf{\#}
&
\textbf{Conservation laws of system}~(\protect\ref{eq:sys1ext_inf})
&
\textbf{Conservation laws of scheme}~(\protect\ref{SamPopScmInfMod1})
&
\textbf{Physics interpretation}
\\
\midrule
1
&
$
\left(\frac{1}{\rho}\right)_t - u_s = 0
$
&
$
\left(\frac{1}{\rho}\right)_t  - \left(u^{(0.5)}\right)_s = 0
$
&
Mass conservation
\\ 
$2^\dag$
&
$
\left(\frac{H^y}{\rho}\right)_t - (H^0 v)_s = 0
$
&
$
\left(\frac{H^y}{\rho}\right)_t - \left(H^0 v\right)_s = 0
$
&
Magnetic flux conservation
\\ 
$3^\dag$
&
$
\left(\frac{H^z}{\rho}\right)_t - (H^0 w)_s = 0
$
&
$
\left(\frac{H^z}{\rho}\right)_t - \left(H^0 w\right)_s = 0
$
&
Magnetic flux conservation
\\ 
4
&
$
u_t + \left(p + \frac{(H^y)^2+ (H^z)^2}{2}\right)_s = 0
$
&
$
u_t + \left(p_-^{(\alpha)} + 
\frac{H^y_-\hat{H}^y_- + H^z_-\hat{H}^z_-}{2} \right)_s = 0
$
&
Momentum conservation
\\ 
$5^\dag$
&
$
v_t - (H^0 H^y)_s = 0
$
&
$
v_t - \left(H^0 H^y_-\right)_{s} = 0
$
&
Momentum conservation
\\ 
$6^\dag$
&
$
w_t - (H^0 H^z)_s = 0
$
&
$
w_t - \left(H^0 H^z_-\right)_{s} = 0
$
&
Momentum conservation
\\ 
7
&
$
\left(t u - x\right)_t + \left[t \left(p + \frac{(H^y)^2+(H^z)^2}{2}\right)\right]_s = 0
$
&
$
\left({\displaystyle \frac{t + \check{t}}{2}} \, u - x\right)_t
+ \left[ t\left(p_-^{(\alpha)} + 
\frac{H^y_-\hat{H}^y_- + H^z_-\hat{H}^z_-}{2} \right)\right]_s = 0
$
&
Center of mass law
\\ 
$8^\dag$
&
$
(t v - y)_t - (t H^0 H^y)_s = 0
$
&
$
\left(t v - y\right)_t
- \left( \hat{t} H^0 \hat{H}^y_- \right)_s = 0
$
&
Center of mass law
\\ 
$9^\dag$
&
$
(t w - z)_t - (t H^0 H^z)_s = 0
$
&
$
\left(t w - z\right)_t
- \left( \hat{t} H^0 \hat{H}^z_- \right)_s = 0
$
&
Center of mass law
\\ 
10
&
$
{\def\arraystretch{2}
\begin{array}{l}
\left(
    \varepsilon
    + \frac{u^2 + v^2 + w^2}{2}
    + \frac{(H^y)^2 + (H^z)^2}{2\rho}
\right)_t +
\\
+ \left[
    \left(p + \frac{(H^y)^2 + (H^z)^2}{2}\right) u
    + E^y H^z - E^z H^y
    \right.
\\
    \left.
    -H^0 (v H^y + w H^z)
\right]_s = 0
\end{array}
} 
$
&
$
{\def\arraystretch{2}
\begin{array}{l}
\left(\check{\varepsilon} + \frac{\check{u}^2 + \check{u}_+^2 + \check{v}^2 + \check{v}_+^2 + \check{w}^2 + \check{w}_+^2}{4}
+ \frac{H^y \check{H}^y + H^z \check{H}^z}{2{\rho}}
\right)_t
\\
+ \left[
    \left( \check{p}_*^{(\alpha)} +
    \frac{(H^y\check{H}^y + H^z\check{H}^z)_*}{2} \right) \check{u}^{(0.5)}
    \right.
    \\
    \left.
    - H^0 (
        \check{v}^{(0.5)} {H}^y_*
        + \check{w}^{(0.5)} {H}^z_*
    )
\right]_s = 0
\end{array}
} 
$
&
Energy conservation
\\
$11^\dag$
&
$
(z v - y w)_t + \left( H^0 (y H^z - z H^y) \right)_s = 0.
$
&
$
\displaystyle
(z v - y w)_t + \left( H^0 (\hat{y} \hat{H}^z_- - \hat{z} \hat{H}^y_-) \right)_s = 0
$
&
Angular momentum conservation
\\
\midrule
\multicolumn{4}{c}{
$
\displaystyle
\varepsilon = \frac{p^{(\alpha)}}{\sum_{k=0}^{\gamma-2} \hat{\rho}^{\gamma-k-1} {\rho}^{k - \gamma + 2}},
\qquad
\gamma \in \mathbb{N}\backslash\{1\}
$
}
\\
\midrule
$12$
&
$
\displaystyle
\left(\frac{p}{\rho^\gamma}\right)_t = 0
$
&
$
\displaystyle
\left(\frac{(\gamma-1)\check{p}^{(\alpha)}}{\sum_{k=0}^{\gamma-2} \rho^{\gamma-k-1} \check{\rho}^{k+1} }\right)_t = 0
$
&
Entropy conservation
\\
\bottomrule
\end{tabularx}
\end{adjustwidth}
\end{table}
\endgroup

\begin{Remark}
\label{rem:clsinf}
Notice that in the case $H^0 = 0$, according to Table~\ref{tab:CLsInf}, scheme~(\ref{SamPopScmInfMod1}) possesses
an infinite set of conservation laws for the following form
\begin{equation}
\displaystyle
\left\{\Phi\left(
\frac{(\gamma-1){p}^{(\alpha)}}{\sum_{k=0}^{\gamma-2} \hat{\rho}^{\gamma-k-1} {\rho}^{k+1}},
\frac{H^y}{\rho}, \frac{H^z}{\rho},
v, w,  y - t v, z - t w
\right) \right\}_t  = 0
\end{equation}
where $\gamma \in \mathbb{N} \backslash \{1\}$ and $\Phi$ is an arbitrary function of its arguments.
\end{Remark}

\begin{Remark}
From (\ref{schemeEntropyCLinf2}) it follows that
\begin{equation} \label{schemeEntropyCLinf2int}
\frac{\check{p}^{(\alpha)}}{\rho\check{\rho}} - S = 0.
\end{equation}
The Taylor series expansion of the latter equation is
\begin{equation}
\frac{p}{\rho^2}  - S + \left[ \frac{p \rho_t}{\rho^3}
    + (\alpha -1)\frac{p_t}{\rho^2} \right]\tau + O(\tau^2) = 0.
\end{equation}
Equation~(\ref{entropyDfferential}) for $\gamma=2$ can be represented as
\begin{equation} \label{entropyDfferentialV2}
\frac{p^{(\alpha)}}{(\rho^{(\alpha)})^2} - \frac{\rho_t^2}{p_t} = 0.
\end{equation}
The corresponding expansion is
\begin{equation}
\frac{p}{\rho^2}
    - \frac{\rho_t^2}{p_t}
    +\left[
        \frac{\rho p_t - 2 \rho_t p}{\rho^3} \, \alpha
        + \frac{(\rho_t p_{tt} - 2 p_t \rho_{tt}) \rho_t}{2 p_t^2}
    \right] \tau
    + O(\tau^2) = 0.
\end{equation}
Equations~(\ref{entropyDfferentialV2})
and~(\ref{schemeEntropyCLinf2int}) approximate the conservation of
entropy with the same order~$O(\tau)$. In contrast
to~(\ref{entropyDfferentialV2}),
approximation~(\ref{schemeEntropyCLinf2int}) can be written in a
divergent form. Thus, it represents a conservation law of the
scheme, while~(\ref{entropyDfferentialV2}) does not. This gives an
advantage in the case of isentropic flows when additional
conservation laws include entropy. Then the expression for the
entropy given by equation~(\ref{schemeEntropyCLinf2int}) can be
considered as a constant and  included into conserved quantities.
Invariant schemes and their conservation laws in case of isentropic
flows are discussed in the next section.
\end{Remark}

\subsubsection{On specific symmetries and conservation laws in case of isentropic flows ($S = \text{const}$)}
\label{sec:isentrop}

According to \cite{bk:DorKozMelKap_PlainFlows_2021}, in case $S = p/\rho^\gamma=\text{const}$,
equations~(\ref{eq:sys1ext_inf}) admit the following symmetries.

\begin{enumerate}

\item
If $H^0 \neq 0$, the admitted Lie algebra is
\begin{equation} \label{isentro_inf_syms_H0_neq_0}
\def\arraystretch{1.75}
\begin{array}{c}
{X_1} = \partial_{t},
\qquad
{X_2} = \partial_{s},
\qquad
{X_3} = \partial_{x},
\qquad
{X_4} = t \partial_{x} + \partial_{u},
\\
{X_5} = z \partial_{y} - y \partial_{z}
    + w \partial_{v} - v \partial_{w}
    + E^z \partial_{E^y} - E^y \partial_{E^z} + H^z \partial_{H^y} - H^y \partial_{H^z},
\\
X_6 = t \partial_{t} + 2 s \partial_{s}
    - u \partial_{u} - v \partial_{v} - w \partial_{w}
    + 2 \rho \partial_{\rho},
\\
X_7 = -s \partial_{s} + x \partial_{x} + y \partial_{y} + z \partial_{z}
+ u \partial_{u} + v \partial_{v} + w \partial_{w}
- 2 \rho \partial_{\rho},
\\
X_8 = q_1(s) \partial_{y},
\qquad
X_9 = q_2(s) \partial_{z},
\qquad
X_{10} = t \partial_{y} + \partial_{v},
\qquad
X_{11} = t \partial_{z} + \partial_{w},
\end{array}
\end{equation}
where $q_1$, $q_2$ are arbitrary functions of~$s$.

Scheme~(\ref{SamPopScmInfMod1}) supplemented by the state equation~(\ref{schemeEntropyCLinfN})
admits all the generators~(\ref{isentro_inf_syms_H0_neq_0}).

\item
In case $H^0=0$, the admitted Lie algebra is
\begin{equation} \label{isentro_inf_syms_H0_eq_0}
\def\arraystretch{1.75}
\begin{array}{c}
{X_1} = \partial_{t},
\quad
{X_2} = \partial_{s},
\quad
{X_3} = \partial_{x},
\quad
{X_4} = t \partial_{x} + \partial_{u},
\quad
X_5 = q_3(s) (H^z \partial_{H^y} - H^y \partial_{H^z}),
\\
X_6 = t \partial_{t} + 2 s \partial_{s} - u \partial_{u} + 2 \rho\partial_{\rho},
\qquad
X_7 = -s \partial_{s} + x \partial_{x} + u \partial_{u} - 2 \rho\partial_{\rho},
\\
X_8 = 2 s \partial_{s} + 2 \rho\partial_{\rho} + p \partial_p + H^y \partial_{H^y} + H^z \partial_{H^z}.
\end{array}
\end{equation}
In case $\gamma=2$, there are two additional generators, namely
\begin{equation} 
X_9 = q_4(s) \rho (\partial_{H^y} - H^y \partial_{p}),
\qquad
X_{10} = q_5(s) \rho (\partial_{H^z} - H^z \partial_{p}).
\end{equation}
Here $q_3$, $q_4$ and $q_5$ are arbitrary functions of~$s$.

Scheme~(\ref{SamPopScmInfMod1}),~(\ref{schemeEntropyCLinfN})
admits all the generators~(\ref{isentro_inf_syms_H0_eq_0}).
However, the scheme does not admit the generators~$X_9$ and~$X_{10}$.
\end{enumerate}

There are the following additional conservation laws for system~(\ref{eq:sys1ext_inf}).

\begin{enumerate}[label=\alph*)]

\item

In case $H^0 \neq 0$, there is an additional conservation law
which corresponds to the generator $\partial_s$

\begin{equation}
\label{SconstCLadd1}
\left(
  \frac{u}{\rho} + \frac{v H^y + w H^z}{H^0 \rho}
\right)_t
+
\left(
    \frac{\gamma S}{\gamma - 1} \rho^{\gamma - 1}
    - \frac{u^2 + v^2 + w^2}{2}
\right)_s = 0.
\end{equation}

\item

Case $H^0 = 0$.

\begin{itemize}

    \item

    The conservation law corresponding to the generator $\partial_s$ is
    \begin{equation}
    \label{SconstCLadd2}
    \left(
       \frac{u}{\rho}
    \right)_t
    +
    \left(
        \frac{\gamma  S_0}{\gamma - 1} \rho^{\gamma - 1}
        - \frac{u^2}{2}
        + \rho B_0
    \right)_s = 0
    \end{equation}
    provided
    \begin{equation}
    \label{CLcondB0}
        S_0 = S = \text{const}
        \qquad
        \text{and}
        \qquad
        B_0 = \frac{(H^y)^2 + (H^z)^2}{\rho^2} = \text{const}.
    \end{equation}
    The latter follows from system (\ref{eq:sys1ext_inf}).
    When conductivity of the medium tends to infinity, the phenomenon of frozen-in magnetic field is observed (see, e.g.~\cite{bk:KulikovskyLubimovMHD1965}). In this case, in the absence of the longitudinal component $H^0$ of the magnetic field, the quantity~$B_0$ which is proportional to the magnetic pressure turns out to preserve along the pathlines.

    \item

    In case $\gamma= 2$, the admitted generator
    \begin{equation}
    \label{addGen3}
        5 X_6 + 3 X_7 - 4 X_8
    \end{equation}
    corresponds to the conservation law
    \begin{equation}
    \label{SconstCLadd3}
    \left(
        5 t \rho A_0
        + 5 t \frac{u^2}{2}
        - s \frac{u}{\rho}
        - 3 x u
    \right)_t
    +
    \left(
        ((5 t u - 3 x) \rho^2 - 2 s \rho) A_0
       + s\frac{u^2}{2}
    \right)_s = 0
    \end{equation}
    provided
    \begin{equation}
    \label{CLcondB1}
    A_0 = \left[
        \frac{S}{\gamma - 1}
        + \frac{(H^y)^2 + (H^z)^2}{2\rho^2}
    \right]_{\gamma=2}
        = \frac{p}{\rho^2} + \frac{(H^y)^2 + (H^z)^2}{2\rho ^2}
        = \text{const}
    \end{equation}
    which follows from system (\ref{eq:sys1ext_inf}).
    \begin{Remark}
    Conservation law~(\ref{SconstCLadd3}) is a basis one. Its partial derivative
    with respect to~$s$ is equivalent to~(\ref{SconstCLadd2}),
    and its partial derivative with respect to~$t$ is
    \begin{equation}
    A_0 (\rho_t + \rho^2 u_s) + u (u_t + (A_0 \rho^2)_s) = 0
    \end{equation}
    which is a combination of (\ref{MassEq}) and (\ref{VelosityEq}) provided~(\ref{CLcondB1}).
    \end{Remark}

\end{itemize}

\end{enumerate}

\bigskip
\par
\bigskip

By virtue of the content of Remark~\ref{rem:clsinf},
one can verify that the finite-difference analogues of~(\ref{CLcondB0}) and (\ref{CLcondB1})
hold along the pathlines for scheme~(\ref{SamPopScmInfMod1}), namely
\begin{equation}
\left(\frac{(H^y)^2 + (H^z)^2}{\rho^2}\right)_t = 0
\quad \text{or} \quad
\left(\frac{H^y\check{H}^y + H^z \check{H}^z}{\rho \check{\rho}}\right)_t = 0
\qquad
\text{if} \quad H^0= 0,
\end{equation}
and
\begin{equation}
\left(
    \frac{p^{(\alpha)}}{\check{\rho}\rho}
    + \frac{(H^y)^2 + (H^z)^2}{2\rho ^2}
\right)_t = 0
\quad \text{or} \quad
\left(
    \frac{p^{(\alpha)}}{\check{\rho}\rho}
    + \frac{H^y\check{H}^y + H^z \check{H}^z}{2\rho \check{\rho}}\right)_t = 0
\qquad
\text{if} \quad H^0 = 0, \; \gamma = 2.
\end{equation}

Scheme~(\ref{SamPopScmInfMod1}) also admits the generators~$\partial_{s}$ and~(\ref{addGen3})
under the same conditions as for the differential case.

\bigskip
\par
\bigskip

Analyzing scheme~(\ref{SamPopScmInfMod1}), one can conclude that for the additional conservation laws~(\ref{SconstCLadd1}), (\ref{SconstCLadd2}) and (\ref{SconstCLadd3}) there is no approximations in terms of rational expressions.
This means that construction of finite difference analogues of the mentioned conservation laws is extremely hard.

\bigskip

Further we restrict ourselves to the case $\gamma=2$ and $S=S_1=\text{const}$, and consider another invariant scheme on an extended finite difference stencil.

We introduce the pressure for the polytropic gas as
\begin{equation}
p = S_1\hat{\rho}\hat{\rho}_+.
\end{equation}
Then, the conservation law of entropy
\begin{equation}
\left(\frac{\check{p}}{\rho\rho_+}\right)_t = (S_1)_t = 0
\end{equation}
is defined by the following invariant expression
\begin{equation}
\frac{p}{\hat{\rho}\hat{\rho}_+}
    = \frac{\check{p}}{\rho\rho_+}
    = \frac{\check{p}_-}{\rho\rho_- }
    = \frac{p_-}{\hat{\rho}\hat{\rho}_-} = S_1.
\end{equation}
The scheme under consideration is based on scheme~(\ref{SamPopScmInfMod1}) and it has the following form
\begin{subequations} \label{SamPopScmInfMod2}
\begin{equation}
\left(\frac{1}{\rho}\right)_t = u^*_s,
\end{equation}
\begin{equation}
\displaystyle
u_t = -\left( p
    + \frac{\hat{H}^y_+ \hat{H}^y + \hat{H}^z_+ \hat{H}^z} {2}\right)_{\bar{s}},
\qquad
v_t = H^0 \hat{H}^y_{s},
\qquad
w_t = H^0 \hat{H}^z_{s},
\end{equation}
\begin{equation} \label{usws}
\left(\frac{H^y}{\rho}\right)_t = H^0 v_{\bar{s}},
\qquad
\left(\frac{H^z}{\rho}\right)_t = H^0 w_{\bar{s}},
\end{equation}
\begin{equation}
\frac{p}{\hat{\rho}\hat{\rho}_+}
    = \frac{\check{p}}{\rho\rho_+}
    = \frac{\check{p}_-}{\rho\rho_- }
    = \frac{p_-}{\hat{\rho}\hat{\rho}_-} = S_1,
\end{equation}
\begin{equation}
x_t = u^*,
\qquad
y_t^+ = v,
\qquad
z_t^+ = w,
\qquad
x_s = \frac{1}{\rho}.
\end{equation}
\end{subequations}
One can verify that the latter scheme is invariant. It admits the same symmetries as scheme~(\ref{SamPopScmInfMod1}),~(\ref{schemeEntropyCLinfN}).
The following quantities hold for~(\ref{SamPopScmInfMod2})
\begin{equation}
\left(
    \frac{H^y_+ H^y + H^z_+ H^z}{\rho\rho_+}
\right)_t  = (B_1)_t = 0
\qquad
\text{if} \quad H^0= 0,
\end{equation}
\begin{equation}
\left(
    \frac{\check{p}}{\rho\rho_+}
    + \frac{H^y_+ H^y + H^z_+ H^z}{2\rho\rho_+}
\right)_t = \frac{1}{2}(2 S_1 + B_1)_t = 0
\qquad
\text{if} \quad H^0 = 0, \; \gamma = 2.
\end{equation}

Scheme (\ref{SamPopScmInfMod2}) possesses the difference analogues of~(\ref{SconstCLadd1}) and~(\ref{SconstCLadd2}), namely
\begin{equation}
\left(
    \frac{u}{\rho}
    + \frac{v_* \hat{H}^y + w_* \hat{H}^z}{H^0 \hat{\rho}}
\right)_t
+ \left(
    2 \hat{\rho}_* S_1
    -\frac{u u_- + \hat{v}^2_- + \hat{w}^2_-}{2}
\right)_s
= 0,
\end{equation}
\begin{equation}
\left(
    \frac{u}{\rho}
\right)_t
+ \left(
    \hat{\rho}_* (2 S_1 + B_1)
    - \frac{u u_-}{2}
\right)_s
= 0.
\end{equation}
To construct the latter conservation law one should use the following relation
\begin{equation}
\frac{\hat{H}^y_+ \hat{H}^y + \hat{H}^z_+ \hat{H}^z} {2}
= \frac{\hat{H}^y_+ \hat{H}^y + \hat{H}^z_+ \hat{H}^z} {2 \hat{\rho}\hat{\rho}_+}
= \frac{1}{2} \hat{\rho}\hat{\rho}_+ B_1.
\end{equation}

\begin{Remark}
The angular momentum and center of mass conservation laws are
\begin{equation}
(z v_- - y w_-)_t + \left( H^0 (\hat{y} \hat{H}^z_- - \hat{z} \hat{H}^y_-) \right)_s = 0,
\end{equation}
\begin{equation}
(t u^* - x)_t + \left(\hat{t} \left[ p_- + \frac{\hat{H}^y_- \hat{H}^y + \hat{H}^z_- \hat{H}^z} {2} \right]_*\right)_s = 0,
\end{equation}
\begin{equation}
(t v - y_+)_t - \left(\hat{t} H^0 \hat{H}^y\right)_s = 0,
\end{equation}
\begin{equation}
(t u - z_+)_t - \left(\hat{t} H^0 \hat{H}^z\right)_s = 0.
\end{equation}
The remaining conservation laws of mass, momentum, magnetic flux and entropy
follow directly from the scheme since it is written in a divergent form.
\end{Remark}

\section{Numerical experiments}
\label{sec:experiments}

In this section, we consider the problem of deceleration of a plasma bunch in a crossed electromagnetic field under the presence and absence of a longitudinal component of magnetic field. We use scheme~(\ref{eq:scheme0ext}), and consider how the conservation laws hold on the solutions of this scheme. In addition to the transverse component~$H^y$ of the magnetic field, we also consider the case of the presence of a longitudinal magnetic field~$H^0\neq 0$.

A plasma bunch is considered, which moves from left to right in a railgun channel. The channel is filled with a relatively cold weakly conducting gas. With the help of an external electric circuit, a strong transverse magnetic field is generated in the channel, which causes the bunch to decelerate. During its motion, the plasma bunch closes
the electric circuit and moves along the background gas; therefore, the magnetic field and pressure at the left boundary of the computational domain are considered equal to zero.
The differential boundary conditions are as follows
\begin{subequations}
\begin{equation}
p(0, t) = 0,
\qquad
H^y(0, t)= 0,
\qquad
H^z(0,t) = 0,
\end{equation}
\begin{equation}
u(S, t) = 0,
\qquad
E^y(S,t) = 0,
\qquad
H^y(S,t) = 4\pi J(t),
\end{equation}
\begin{equation} \label{BCsDiff1}
\displaystyle
L_0 \frac{dJ}{dt} + R_0 J - V(t) + E^z(S,t) = 0,
\end{equation}
\begin{equation} \label{BCsDiff2}
\frac{dV}{dt} = -J/C_0, \quad V(0) = V_0, \quad J(0)=0,
\end{equation}
\end{subequations}
where $0 \leqslant s \leqslant S$ and $0 \leqslant t \leqslant t_{\max}$, $S$ is the total mass of the gas,
$J$ and $V$ are current and voltage, and~$C_0$,~$L_0$,~$R_0$ are the external circuit parameters.
The boundary conditions~(\ref{BCsDiff1}) and~(\ref{BCsDiff2}) are approximated in the same way as in~\cite{bk:SamarskyPopov_book[1992]}, namely
\begin{equation}
L_0 J_t + R_0 J^{(0.5)} - V^{(0.5)} + {E^z_M}^{(0.5)} = 0,
\qquad
V_t = -J^{(0.5)}/C_0,
\end{equation}
where $M=\lfloor{S/h}\rfloor$.

All calculations were carried out using the dimensionless version of scheme~(\ref{eq:scheme0ext}) with the value of the coefficient~$\kappa=4\pi$.
For the dimensionless form of the scheme, the initial conditions are: $\rho_0=1.0$, $p_0=0.0056$, $R_0=1.17$, $C_0=1.64$, $L_0=0.0035$, $S=4.0$, the temperature of the plasma~$T_0=3.0$, the initial speed of the plasma bunch~$u_0 = 0.75$. The gas is considered polytropic with~$\gamma=5/3$. The uniform mesh steps are~$h=0.067$ and~$\tau=0.003$, and~$t_{\max} = 0.7$.
The initial voltage~$V_0$ is varied between~$1.67$ and~$2.6$ which approximately correspond to the voltage~$650$ and $1000$~V. In experiments where the longitudinal magnetic field~$H^0$ is present, a value close to~1 is taken for~$H^0$. In the calculations, a linear artificial viscosity is used, with a viscosity coefficient~$\nu = 2 h$.

The problem under consideration is close to the problem described in~\cite{bk:TLayersPreprintMain1973} in which,
however, tabulated real plasma parameters, including electrical conductivity, were used. In our problem, we used the ideal gas equation and an exponential conductivity function.

Scheme~(\ref{eq:scheme0ext}) is implemented using the iterative methods described in~\cite{bk:SamarskyPopov_book[1992]}. In this case, the scheme equations are divided into two parts, dynamic and magnetic. The dynamic part is preliminarily linearized using the Newton method, and for the magnetic part a flow version of the sweep method is used~\cite{bk:DEGTYAREV1968252},
which is well suited for the case of finite conductivity, especially when its values are small.
The bunch motion is modeled by a shock wave.
Conductivity~$\sigma$
of the plasma bunch is proportional to~$T^{3/2}$, and the conductivity function~$\sigma=\sigma(\rho,T)$ is very sensitive to the density~$\rho$ in such a way that in the rarefied background gas region it has values close to zero.

\medskip

Three essentially different cases are considered:
\begin{enumerate}
  \item
  The bunch is decelerated using a transverse magnetic field~$H^y$ at a relatively low voltage in the circuit.

  \item
  The bunch is decelerated using a transverse magnetic field at a high voltage in the circuit.

  \item
  A rather strong longitudinal magnetic field~$H^0$ is added to the previous case. (Calculations show that a weak longitudinal magnetic field has little effect on the experimental results.)
\end{enumerate}
In all cases, at the initial moment of time, the gas particles are given a small constant transverse velocity~$v>0$.
This is necessary in order to track the influence of the longitudinal magnetic field on the transverse component of the particle velocity, which should be observed only in the third numerical experiment.

\begin{figure}[H]
\begin{adjustwidth}{-\extralength}{0cm}
\includegraphics[width=0.75\linewidth]{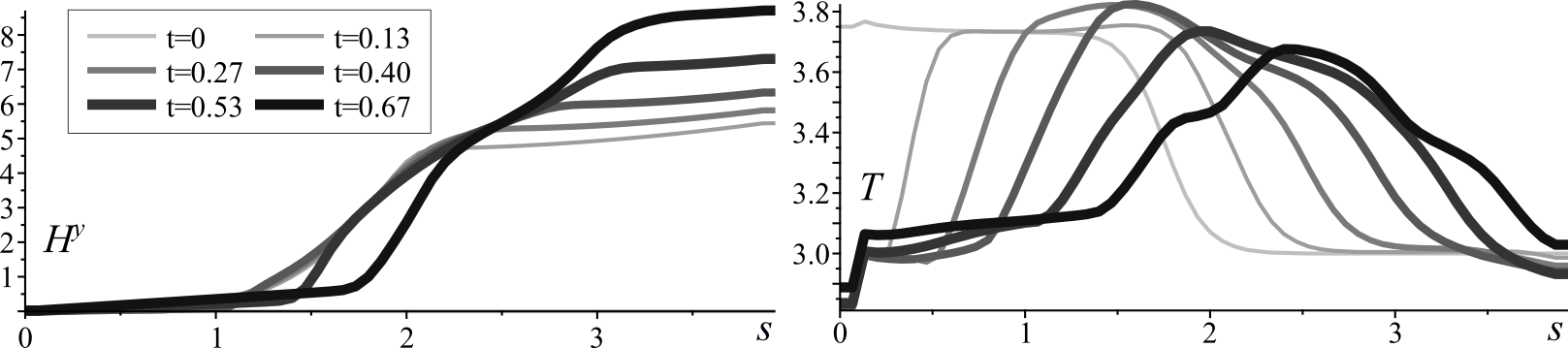}
\centering
\end{adjustwidth}
\caption{Evolution of the magnetic field~$H^y$ and the temperature~$T$ for the first experiment.}
\label{fig:case1}
\end{figure}

\begin{figure}[H]
\begin{adjustwidth}{-\extralength}{0cm}
\centering
\includegraphics[width=0.75\linewidth]{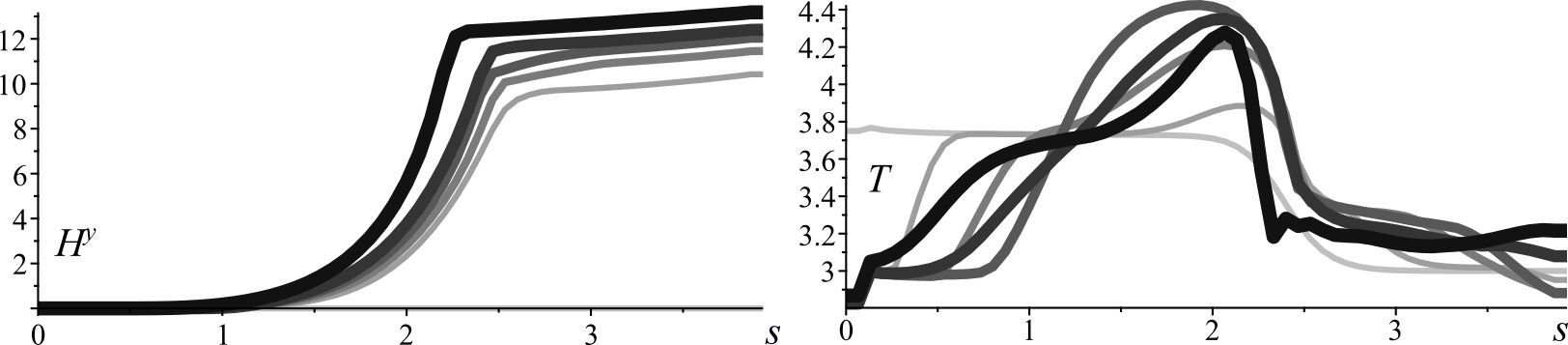}
\end{adjustwidth}
\caption{Evolution of the magnetic field~$H^y$ and the temperature~$T$ for the second experiment.}
\label{fig:case2}
\end{figure}

\begin{figure}[H]
\begin{adjustwidth}{-\extralength}{0cm}
\centering
\includegraphics[width=0.75\linewidth]{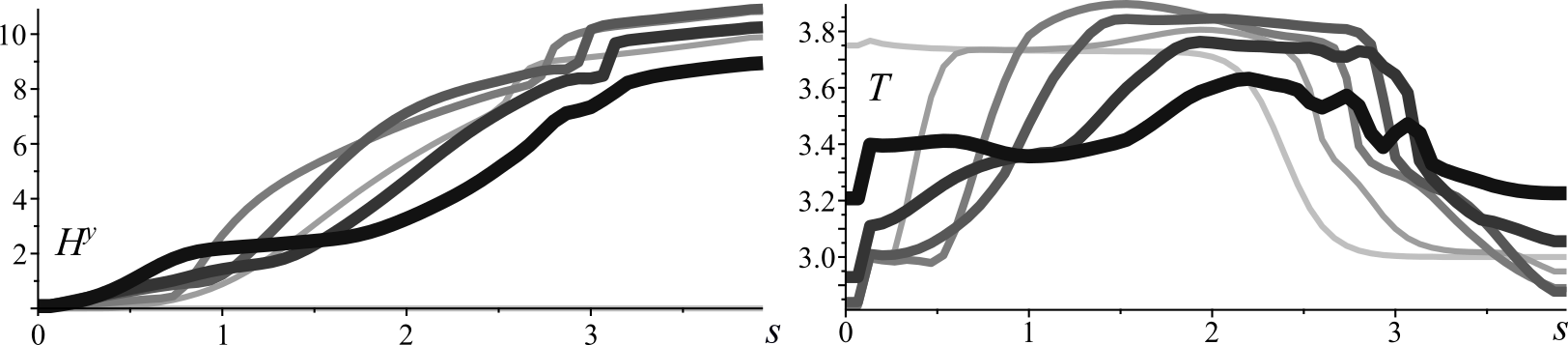}
\end{adjustwidth}
\caption{Evolution of the magnetic field~$H^y$ and the temperature~$T$ for the third experiment.}
\label{fig:case3}
\end{figure}

Figure~\ref{fig:case1} shows the evolution of the magnetic field and plasma temperature in the first experiment. The magnetic field is not strong enough to stop the bunch. If the bunch reaches the right boundary of the computational domain, the reflection of the wave can be observed due to the boundary condition $u(S,t) = 0$.
Figure~\ref{fig:case2} shows the second case where the transverse magnetic field is strong enough. The plasma bunch is decelerated by the magnetic field and after a short period of time begins to move backward. Adding a sufficiently strong longitudinal magnetic field~$H^0$ to the previous experiment leads to an intermediate picture: the magnetic field is ``smeared'' over the computational domain, the plasma deceleration process is not as intense as in the previous case, and is inhomogeneous along the mass coordinate, which leads to a kind of fragmentation of the temperature profile~(see Figure~\ref{fig:case3}).


\begin{figure}[H]
\begin{adjustwidth}{-\extralength}{0cm}
\centering
\includegraphics[width=0.8\linewidth]{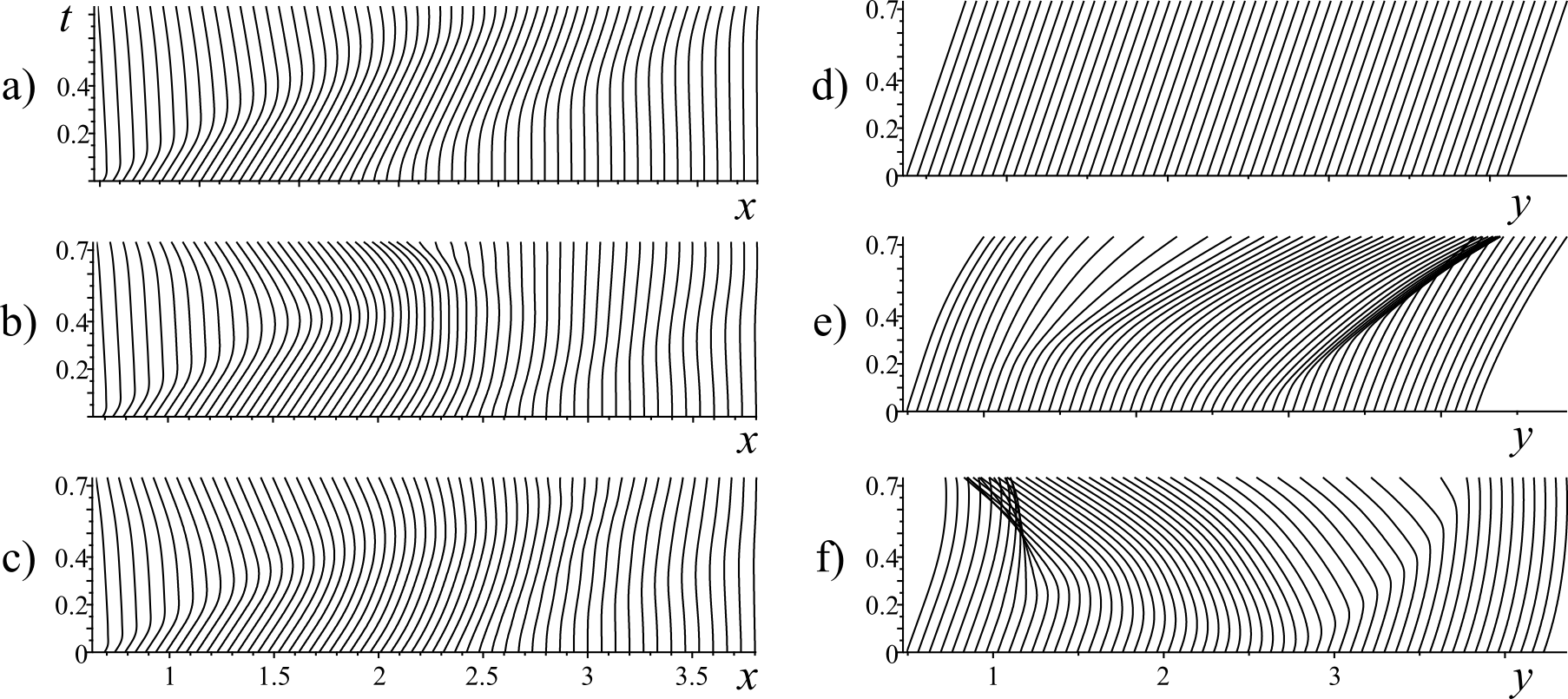}
\end{adjustwidth}
\caption{a), b) and c) show $x$-trajectories for the experiments 1, 2 and 3; d) $y$-trajectories for~$H^0=0$; e) $y$-trajectories for~$H^0>0$; f) $y$-trajectories for~$H^0<0$.}
\label{fig:traj}
\end{figure}

In Figure~\ref{fig:traj} the trajectories of particles under the action of magnetic fields are shown. The left part~(Figure~\ref{fig:traj},~a)--c)) shows~$x$-trajectories of particles for three experiments. The right side of the figure shows~$y$-trajectories associated with the transverse velocity component~$v$. Figure~\ref{fig:traj},~d) corresponds to the first and second experiments where~$v$ has a constant value and~$H^0=0$. Figure~\ref{fig:traj},~e) and Figure~\ref{fig:traj},~f) correspond to the third experiment at~$H^0>0$ and~$H^0<0$ where under the action of the longitudinal magnetic field the transverse velocity component increases or slows down accordingly. Notice that the choice of sign of the value~$H^0$ otherwise does not affect the results of the third experiment.

\begin{figure}[H]
\begin{adjustwidth}{-\extralength}{0cm}
\centering
\includegraphics[width=0.5\linewidth]{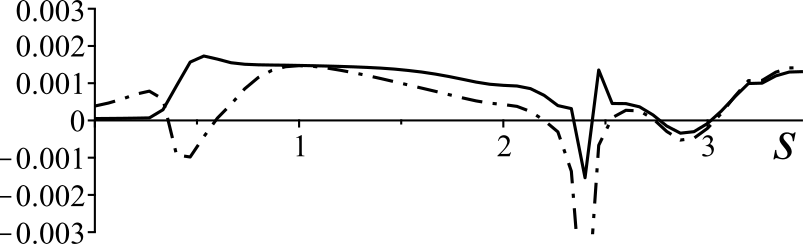}
\end{adjustwidth}
\caption{Conservation laws of energy (solid line) and $y$-flux (dashed line) at $t=0.64$.}
\label{fig:cls1}
\end{figure}

\begin{figure}[H]
\begin{adjustwidth}{-\extralength}{0cm}
\centering
\includegraphics[width=0.5\linewidth]{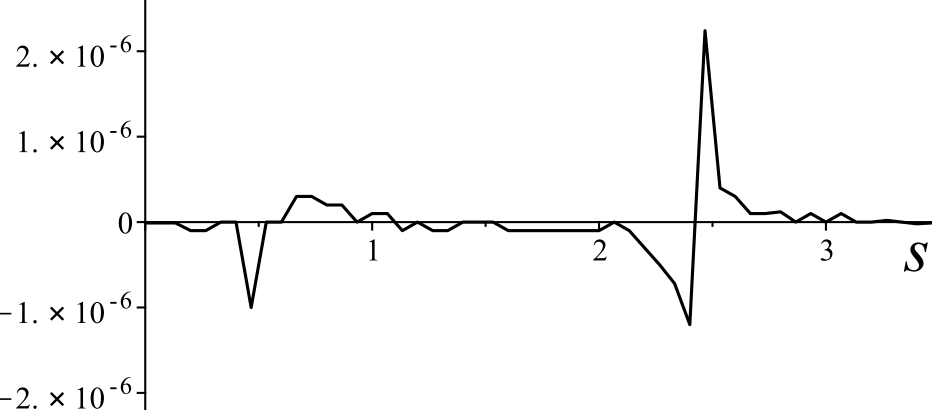}
\end{adjustwidth}
\caption{Conservation law of $x$-momentum at $t=0.64$.}
\label{fig:cls2}
\end{figure}

\begin{figure}[H]
\begin{adjustwidth}{-\extralength}{0cm}
\centering
\includegraphics[width=0.5\linewidth]{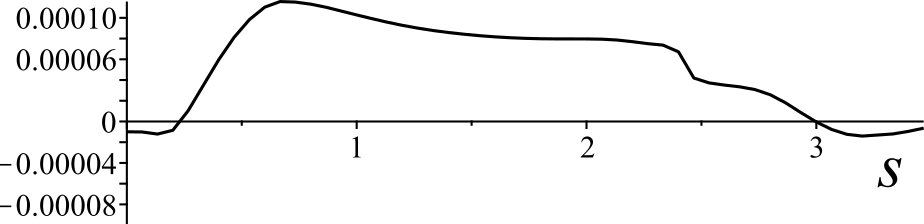}
\end{adjustwidth}
\caption{Conservation law of center-of-mass along axis~$x$ at $t=0.64$.}
\label{fig:cls3}
\end{figure}

In Figures~\ref{fig:cls1}--\ref{fig:cls3} the finite-difference conservation laws of energy, magnetic flux (along the~$y$ axis), momentum and center-of-mass motion (along the~$x$ axis) are given for the selected moment of time, when the interaction of magnetic fields and the plasma bunch is already quite intense.
The results are provided only for the third experiment, since in other cases the control of conservation laws gives similar results. The accurate enough preservation of the conservation laws on solutions is due to the conservativeness of scheme~(\ref{eq:scheme0ext}).


\section{Conclusion}

Finite-difference schemes for MHD equations in the case of plane one-dimensional flows are considered.
The Samarsky--Popov classical scheme for the case of finite conductivity is taken as a starting point.
Symmetries and conservation laws of this scheme are investigated.
It is shown that the scheme admits the same symmetries as
the original differential model.
It also has difference analogues of the conservation laws of the original model.
In addition to the conservation laws previously known for the scheme,
new conservation laws are given, which are obtained on the basis of the group classification recently carried out in~\cite{bk:DorKozMelKap_PlainFlows_2021}.

The classical Samarskiy--Popov scheme is generalized to the case of
arbitrary vectors of electric and magnetic fields, as well as to the
case of infinite conductivity.
In the case of finite conductivity
the scheme possesses difference analogues of all differential local
conservation laws obtained
in~\cite{bk:DorKozMelKap_PlainFlows_2021}, some of which were not
previously known. In the case of infinite conductivity,
straightforward generalization of the scheme leads to a scheme that
does not preserve angular momentum. The proposed modification
makes it possible to obtain an invariant scheme that also possesses the conservation law of angular momentum.
In addition, it is shown how to approximate the
equation of state for a polytropic gas to preserve the entropy along
the pathlines on the extended stencil for two time layers.

A numerical implementation of the generalized Samarsky-Popov scheme for the case of finite conductivity is performed for the problem of deceleration of a plasma bunch by crossed electromagnetic fields. Various cases of the action of fields on a plasma are considered. Calculations show that the finite-difference conservation laws are preserved on the solutions of the scheme quite accurately.

\vspace{6pt}

\authorcontributions{Conceptualization, V.D.; methodology, V.D. and E.K.; software, E.K.; validation, E.K; investigation, V.D. and E.K.; data curation, E.K.; writing---original draft preparation, V.D. and E.K.; writing---review and editing, V.D. and E.K.; visualization, E.K.; supervision, V.D.; project administration, V.D.; funding acquisition, V.D. All authors have read and agreed to the published version of the manuscript.}

\funding{This research was supported by Russian Science Foundation Grant No.~18-11-00238
``Hydrodynamics-type equations: symmetries, conservation laws, invariant difference
schemes''.}

\dataavailability{The data presented in this study are available on request from the corresponding author.}

\acknowledgments{The authors thank E.~Schulz and S.V.~Meleshko for valuable discussions.
E.K. sincerely appreciates the hospitality of the Suranaree University of Technology.}

\conflictsofinterest{The authors declare no conflict of interest.}

\reftitle{References}

\begin{adjustwidth}{-\extralength}{0cm}

\end{adjustwidth}

\end{document}